\documentclass[a4paper]{article}

\usepackage{epsfig}
\usepackage[francais,english]{babel}
\usepackage[latin1]{inputenc}
\usepackage[tbtags]{amsmath}
\usepackage{amsthm}
\usepackage{amssymb}
\usepackage{pifont}
\usepackage[all]{xy}

\usepackage{vmargin}

\newtheorem{thintro}{Theorem}
\newtheorem{theo}{Theorem}[subsection]
\newtheorem{lemma}[theo]{Lemma}
\newtheorem{prop}[theo]{Proposition}
\newtheorem{cor}[theo]{Corollary}
\theoremstyle{definition}
\newtheorem{deftn}[theo]{Definition}
\newtheorem{question}{Question}

\def\leq{\leqslant}
\def\geq{\geqslant}

\def\N{\mathbb{N}}
\def\Z{\mathbb{Z}}
\def\Q{\mathbb{Q}}

\def\F{\mathbb{F}}
\def\O{\mathcal{O}}
\def\E{\mathcal{E}}
\def\pa#1{\left(#1\right)}

\def\acco#1{\left\{ #1 \right\}}
\def\brac#1{\left< #1 \right>}

\def\epsilon{\varepsilon}

\def\calC{\mathcal{C}}
\def\calF{\mathcal{F}}
\def\calG{\mathcal{G}}
\def\calL{\mathcal{L}}
\def\calM{\mathcal{M}}
\def\calS{\mathcal{S}}

\def\frakM{\mathfrak{M}}
\def\frakN{\mathfrak{N}}
\def\calK{\mathcal{K}}
\def\frakK{\mathfrak{K}}
\def\frakC{\mathfrak{C}}
\def\frakI{\mathfrak{I}}

\def\si{\mathfrak{S}}

\def\id{\text{\rm id}}
\def\im{\text{\rm im}\:}
\def\coim{\text{\rm coim}\:}
\def\ker{\text{\rm ker}\:}
\def\coker{\text{\rm coker}\:}
\def\hom{\text{\rm Hom}}
\def\Frac{\text{\rm Frac}\:}
\def\card{\text{\rm Card}\:}
\def\Fil{\text{\rm Fil}}
\def\ur{\text{\rm ur}}
\def\sep{\text{\rm sep}}

\makeatletter
\def\virgule#1{%
   \edef\@temp{#1}%
   \ifx\@temp\empty{}\else{#1,}\fi}%
\makeatother

\def\Rep{\text{\rm Rep}}
\def\Repg{\Rep_\infty^{[0,1]}}
\def\pModphiN#1#2{\text{\rm 'Mod}^{\virgule{#1}\phi,N}_{/#2}}
\def\ModphiN#1#2{\text{\rm Mod}^{\virgule{#1}\phi,N}_{/#2}}
\def\Modphin#1#2{\text{\rm Mod}^{\virgule{#1}\phi,(N)}_{/#2}}
\def\pModphi#1#2{\text{\rm 'Mod}^{\virgule{#1}\phi}_{/#2}}
\def\Modphi#1#2{\text{\rm Mod}^{\virgule{#1}\phi}_{/#2}}
\def\Maxphi#1#2{\text{\rm Max}^{\virgule{#1}\phi}_{/#2}}
\def\MaxphiN#1#2{\text{\rm Max}^{\virgule{#1}\phi,N}_{/#2}}
\def\Minphi#1#2{\text{\rm Min}^{\virgule{#1}\phi}_{/#2}}
\def\pM{\text{\rm '}M}
\def\pT{\text{\rm '}T}

\def\st{\text{\rm st}}
\def\cris{\text{\rm cris}}
\def\qcris{\text{\rm qst}}
\def\Ast{\hat A_\st}
\def\Acris{A_\cris}

\def\Max{\text{\rm Max}}
\def\max{\text{\rm max}}
\def\imax{\iota_\max}
\def\Min{\text{\rm Min}}
\def\min{\text{\rm min}}
\def\imin{\iota_\min}

\title{Quasi-semi-stable representations}
\author{Xavier Caruso and Tong Liu}

\begin{document}

\maketitle

\begin{abstract}
Fix $K$ a $p$-adic field and denote by $G_K$ its absolute Galois group.
Let $K_\infty$ be the extension of $K$ obtained by adding $p^n$-th
roots of a fixed uniformizer, and $G_\infty \subset G_K$ its absolute
Galois group.
In this article, we define a class of $p$-adic torsion representations
of $G_\infty$, named \emph{quasi-semi-stable}. We prove that these
representations are ``explicitly'' described by a certain category of
linear algebra objects.
The results of this note should be consider as a first step in the
understanding of the structure of quotients of two lattices in a
crystalline (resp. semi-stable) Galois representation.
\end{abstract}

\tableofcontents

\section*{Introduction}

Let $p$ be an odd prime number and $k$ a perfect field of characteristic $p$.
Put $W = W(k)$ the ring of Witt vectors with coefficients in $k$, and
$K_0 = \Frac W$. Denote by $\sigma$ the Frobenius on $k$, $W$ and $K_0$.
Let $K$ be a totally ramified extension of $K$ of degree $e$ and $\O_K$
its ring of integers. Fix $\pi$ an uniformizer of $\O_K$. We denote by
$\bar K$ an algebraic closure of $K$, by $\O_{\bar K}$ its ring of
integers and by $G_K$ its absolute Galois group. Fix a sequence
$(\pi_n)$ of elements of $\bar K$ satisfying $\pi_0 = \pi$ and
$\pi_{n+1}^p = \pi_n$. Put $K_n = K(\pi_n)$, $K_\infty = \bigcup_{n \in
\N} K_n$ and denote by $G_\infty \subset G_K$ the absolute Galois group
of $K_\infty$.

We would like to study representations that can be written as a quotient
of two lattices in a crystalline or semi-stable representation. For
technical reason we have to make an assumption on Hodge-Tate weights,
that is they all belong to $\{0, \ldots, r\}$ for an integer $r < p-1$.
The theory of Breuil modules then gives a description of these lattices
in term of linear algebra: there exists a category $\ModphiN r S$ that
is dually equivalent to those whose objects are these lattices. By
mimicing the definition of
$\ModphiN r S$, one can construct a category of torsion objects $\ModphiN r
{S_\infty}$ equipped with a contravariant functor $T_\st$ with values in
the category of Galois representations. When $er < p-1$, we can prove
that $\ModphiN r {S_\infty}$ is an abelian category and $T_\st$ is fully
faithful (see \cite{caruso-crelle}). However, these assertions become
false if the assumption $er < p-1$ is removed. In this article, we draw
a picture of the structure of all this stuff in a slighty different
situation. Precisely, we remove the operator $N$ (that appears in the
subscript $\ModphiN r S$) and study a new category so-called $\Modphi
r S$. It is endowed with a functor $T_\qcris$ with values in a certain
category of $G_\infty$-representations, that we call \emph{quasi-semi-stable}.
The following theorem gathers many important results of structure proved 
in this paper.

\begin{thintro}
\label{theo:overS}
Let $\calM \in \Modphi r {S_\infty}$. There exists a unique couple
$(\Max^r(\calM), \imax^\calM : \calM \to \Max^r(\calM))$ (resp.
$(\Min^r(\calM), \imin^\calM : \Min^r (\calM) \to \calM)$) (where
$\Max^r(\calM)$, $\Min^r(\calM)$ are objects of $\Modphi r {S_\infty}$ and
$\imax^\calM$, $\imin^\calM$ morphisms in this category) such that:
\begin{itemize}
\item the morphism $T_\qcris(\imax^\calM)$ (resp. $T_\qcris(\imin^\calM)$)
is an isomorphism;
\item for any $\calM' \in \Modphi r {S_\infty}$ endowed with a morphism $f
: \calM \to \calM'$ (resp. $f : \calM' \to \calM$) such that
$T_\qcris(f)$ is an isomorphism, there exists a unique $g : \calM' \to
\Max^r(\calM)$ (resp. $g : \Min^r(\calM) \to \calM'$) such that $g \circ f =
\imax^\calM$ (resp. $f \circ g = \imin^\calM$).
\end{itemize}
This property gives rise to a functor $\Max^r : \Modphi r {S_\infty} \to
\Modphi r {S_\infty}$ (resp. $\Min^r : \Modphi r {S_\infty} \to \Modphi r
{S_\infty}$) which satisfies $\Max^r \circ \Max^r = \Max^r$ (resp. $\Min^r
\circ \Min^r = \Min^r$). Its essential image $\Maxphi r {S_\infty}$ (resp.
$\Minphi r {S_\infty}$) is an abelian category. The functor $\Max^r : \Modphi r
{S_\infty} \to \Maxphi r {S_\infty}$ (resp. $\Min^r : \Modphi r {S_\infty} \to
\Minphi r {S_\infty}$) is exact and a left adjoint (resp. a right adjoint)
to the inclusion $\Maxphi r {S_\infty} \to \Modphi r {S_\infty}$. (resp.
$\Minphi r {S_\infty} \to \Modphi r {S_\infty})$. The restriction of
$T_\qcris$ on $\Maxphi r {S_\infty}$ (resp. on $\Minphi r {S_\infty}$) is
fully faithful. Its essential image is stable under quotients and
subobjects. Moreover, the functor $\Max^r : \Modphi r {S_\infty} \to \Maxphi r
{S_\infty}$ (resp. $\Min^r : \Modphi r{S_\infty} \to \Minphi r {S_\infty}$)
realizes the localization of $\Modphi r {S_\infty}$ with respect to
morphisms $f$ such that $T_\qcris(f)$ is an isomorphism.

Furthermore, functors $\Max^r$ and $\Min^r$ induce exact equivalences of
categories between $\Minphi r {S_\infty}$ and $\Maxphi r {S_\infty}$,
quasi-inverse one to the other.
\end{thintro}

If $r=1$, quasi-semi-stable representations are linked with geometry. In
this case, the category $\Modphi r {S_\infty}$ is dually equivalent to the
category of finite flat group schemes over $\O_K$ killed by a power of
$p$ (see \cite{breuil-group}). Under this equivalence, the functor
$\Min^r$ (resp. $\Max^r$) corresponds to the maximal (resp. minimal) models
defined by Raynaud in \cite{raynaud}. The following result is then a
direct consequence of theorem \ref{theo:overS}.

\begin{thintro}
\label{theo:raynaud}
The category of minimal (resp. maximal) finite flat group schemes
over $\O_K$ killed by a power of $p$ is abelian.
\end{thintro}

Finally, always in the case $r=1$, we can derive from our results a
new proof of the following theorem.

\begin{thintro}
\label{theo:breuil}
Let $\calG$ and $\calG'$ two finite flat group schemes over $\O_K$
killed by a power of $p$. Put $T = \calG(\bar K)$, $T' = \calG'(\bar
K)$ and consider $f : T \to T'$ a $G_\infty$-equivariant map. Then $f$
is $G_K$-equivariant.
\end{thintro}

Unfortunately, if $r > 1$, quasi-semi-stable representations do not have
anymore a geometric interpretation. Then, it is difficult to derive
concrete results from theorem \ref{theo:overS} in general. Actually,
theorem \ref{theo:overS} should be seen as a preliminary for the study
of the more interesting category $\ModphiN r {S_\infty}$; a first part
of this work will be achieved in a forthcoming paper (see
\cite{caruso-forthcoming}).

\medskip

Now, we detail the structure of the article. First, we recall
definitions of categories of Breuil modules. This allows us to explain
more precisely and more clearly our motivations and results. In the
second section, we introduce the category $\Modphi r {\si_\infty}$ and
we prove that it is equivalent to the category $\Modphi r {S_\infty}$.
This result is interesting because it will be easier to work with
objects of $\Modphi r {\si_\infty}$. Section \ref{sec:const} is devoted
to the study of the structure of $\Modphi r {\si_\infty} = \Modphi r
{S_\infty}$: essentially we give a proof of theorem \ref{theo:overS}.
Then, we assume $r = 1$ and show how the previous theory easily imply
theorem \ref{theo:breuil}. The paper ends with some perspectives and
open questions.

\section{Motivations and settings}

Since, in the rest of the paper, we will make an intensive use of Breuil
modules, we choose to gather below all basic definitions about it.
Maybe, the reader may skip it in a first time and come back after when 
objects are really used.

\subsection{Breuil modules}
\label{subsec:breuilmod}

\emph{Fix an integer $r < p-1$.} Recall that $\pi$ is a fixed uniformizer.
Denote by $S$ the $p$-adic completion of the PD-envelope of $W[u]$ with
respect to the kernel of the surjection $W[u] \to \O_K$, $u \mapsto \pi$
(and compatible with the canonical divided powers on $pW[u]$). This
ideal is principal generated by $E(u)$, the minimal polynomial of $\pi$
over $K_0$. The ring $S$ is endowed with the canonical filtration
associated to the PD-envelope and with two endomorphisms:
\begin{itemize}
\item a Frobenius $\phi$: it is the unique continuous map
$\sigma$-semi-linear which sends $u$ to $u^p$
\item a monodromy operator $N$: it is the unique continuous map
$W$-linear that satisfies the Leibniz rule and sends $u$ to $-u$.
\end{itemize}
They satisfy $N \phi = p \phi N$. We have
$\phi(\Fil^r S) \subset p^r S$ (recall $r < p-1$) and we define $\phi_r
= \frac \phi{p^r} : \Fil^r S \to S$. Put $c = \phi_1(E(u))$: it is a
unit in $S$.

First, we define a ``big'' category $\pModphiN r S$ whose objects are the
following data:
\begin{enumerate}
\item a $S$-module $\calM$;
\item a submodule $\Fil^r \calM \subset \calM$ such that $\Fil^r S \, \calM
\subset \Fil^r \calM$;
\item a $\phi$-semi-linear map $\phi_r : \Fil^r \calM \to \calM$;
\item a $W$-linear map $N : \calM \to \calM$ such that:
\begin{itemize}
\item (Leibniz condition) $N(sx) = s N(x) + N(s) x$ for all $s\in S$, $x \in
\calM$
\item (Griffiths transversality) $E(u) N(\Fil^r \calM) \subset \Fil^r
\calM$
\item the following diagram is commutative:
$$\xymatrix @C=50pt {
\Fil^r \calM \ar[r]^-{\phi_r} \ar[d]_{E(u)N} & \calM \ar[d]^{cN} \\
\Fil^r \calM \ar[r]^-{\phi_r} & \calM  }$$
\end{itemize}
\end{enumerate}
Morphisms in $\pModphiN r S$ are whose that are $S$-linear and
compatible with $\Fil^r$, $\phi_r$ and $N$. There exists in $\pModphiN r
S$ a notion of exact
sequence: a sequence $0 \to \calM' \to \calM \to \calM'' \to 0$ is said
exact if both sequences $0 \to \calM' \to \calM \to \calM'' \to 0$ and
$0 \to \Fil^r \calM' \to \Fil^r \calM \to \Fil^r \calM'' \to 0$ are
exact as sequences of $S$-modules.

\medskip

Now, we are ready to define full subcategories of $\pModphiN r S$. The
first one is the category of \emph{strongly divisible modules}, namely
$\ModphiN r S$: it consists of objects $\calM \in \pModphiN r S$ satisfying
the following conditions:
\begin{itemize}
\item the module $\calM$ is free of finite rank over $S$;
\item the quotient $\calM / \Fil^r \calM$ have no $p$-torsion;
\item the image of $\phi_r$ generates $\calM$ (as an $S$-module).
\end{itemize}
The second category is $\ModphiN r {S_1}$: these objects are the $\calM
\in \pModphiN r S$ such that
\begin{itemize}
\item the module $\calM$ is free of finite rank over $S_1 = S/pS$;
\item the image of $\phi_r$ generates $\calM$ (as an $S$-module).
\end{itemize}
Finally, let $\ModphiN r {S_\infty}$ be the smallest subcategory
of $\pModphiN r S$ containing $\ModphiN r {S_1}$ and stable under extensions
(\emph{i.e.} if $0 \to \calM' \to \calM \to \calM'' \to 0$ is an exact
sequence in $\pModphiN r S$ and if $\calM'$ and $\calM''$ are objects of
$\ModphiN r {S_\infty}$, then $\calM$ is also).

\medskip

The three former categories are equipped with a contravariant functor
$T_\st$ with values in the category of $\Z_p$-representations of $G_K$.
On $\ModphiN r S$, it is defined by the formula
$$T_\st (\calM) = \hom_{\pModphiN r S} (\calM, \Ast)$$
where $\Ast$ is a certain period ring, object of $\pModphiN r S$ endowed
with an action of $G_K$. We refer to \cite{breuil-annals} (\S 3.1.1) for
the precise definition of $\Ast$. On the category $\ModphiN r {S_\infty}$
it is defined by
$$T_\st (\calM) = \hom_{\pModphiN r S} (\calM, \Ast \otimes_{\Z_p}
\Q_p/\Z_p).$$

We define similarly categories $\pModphi r S$, $\Modphi r S$, $\Modphi r
{S_1}$ and $\Modphi r {S_\infty}$ by forgetting the operator $N$. The
three last ones are equipped with a functor $T_\qcris$ with values in the
category of $\Z_p$-representations of $G_\infty$\footnote{$T_\qcris(\calM)$
is not endowed with an action of $G_K$ since this group does not act
trivially on $u \in \Acris$.} (defined in the introduction): definitions
are obtained by replacing the period ring $\Ast$ by $\Acris$.
We have a collection of forgetting functors, and if $\calM$ is an object
of $\ModphiN r S$ (resp. $\ModphiN r {S_\infty}$), we have a canonical
and functorial $G_\infty$-equivariant isomorphism
\begin{equation}
\label{eq:tsttqst}
T_\st(\calM) \simeq T_\qcris(\calM)
\end{equation}
(see lemma 2.3.1.1 of \cite{breuil-invent}).

\subsection{Aim of the paper}

Semi-stable $\Q_p$-representations of $G_K$ are classified by (weakly)
admissible filtered $(\varphi, N)$-modules (see \cite{fontaine}). Our
motivations are to describe quotients of two lattices in such
representations, in term of linear algebra. If the Hodge-Tate weights of
the semi-stable representations is in $\{0, \ldots, r\}$, such a
description exists for lattices (stable by $G_K$):

\begin{theo}[Liu, \cite{liu}]
The functor $T_\st$ from $\ModphiN r S$ to the category of lattices in
semi-stable representations with Hodge-Tate weights in $\{0, \ldots,
r\}$ is an anti-equivalence.
\end{theo}

Furthermore, we have the following lemma:

\begin{lemma}
\label{lem:quotfortdiv}
Let $\calM' \subset \calM$ be two strongly divisible modules such that
$\calM' \otimes_{\Z_p} \Q_p \simeq \calM \otimes_{\Z_p} \Q_p$ and
$\Fil^r \calM' \otimes_{\Z_p} \Q_p \simeq \Fil^r \calM \otimes_{\Z_p}
\Q_p$. Then $\calM/\calM'$ is an object of $\pModphiN r S$ and the
following sequence of $G_K$-representations:
$$0 \to T_\st(\calM) \to T_\st(\calM') \to \hom_{\pModphiN r S}
(\calM/\calM', \Ast \otimes_{\Z_p} \Q_p/\Z_p) \to 0$$
is exact.
\end{lemma}

\begin{proof}
The argument is the same as in lemma V.4.2.4 of \cite{caruso-thesis}.
\end{proof}

\noindent
Knowing this, we can draw a plan to study our representations:
\begin{enumerate}
\item recognize objects in $\pModphiN r S$ that can be written as a
quotient of two divisible modules as in lemma \ref{lem:quotfortdiv};
\item study the functor $\hom_{\pModphiN r S} (\text{---}, \Ast
\otimes_{\Z_p} \Q_p/\Z_p)$ on this subcategory.
\end{enumerate}

\medskip

The aim of this article is to explain how we can lead to end the
previous plan in the case of objects of $\pModphi r S$ and $\Modphi r S$
(\emph{i.e.} without $N$). Precisely we prove that the category of
torsion quotients of two objects of $\Modphi r S$ is exactly the category
$\Modphi r {S_\infty}$, and then theorem \ref{theo:overS}.

\medskip

We can imagine that a representation arising from an object of $\Modphi r
S$ should be just a lattice in a crystalline representation, but
unfortunately the situation is quite more complicated. Lattices in
crystalline representations correspond to objects of $\ModphiN r S$ for
which $N(\calM) \subset (uS + \Fil^1 S) \calM$. Let's call $\Modphin r S$
their subcategory. We can see easily that a $N$ satisfying the previous
condition is necessary unique. However, the following lemma shows that
it does not exist in general.

\begin{lemma}
Assume $r \geq 2$ and consider $\calM$ the object of $\Modphi r S$
defined by the following equations :
\begin{enumerate}
\item $\calM = S e_1 \oplus S e_2$ ;
\item $\Fil^r \calM = E(u)^{r-2} e_1 S + E(u)^r e_2 S + \Fil^p S \: 
\calM$ ;
\item $\phi(e_1) = p^2(e_1 + u e_2)$ and $\phi(e_2) = u e_1 + e_2$.
\end{enumerate}
Then, it is impossible to equip $\calM$ with a monodromy operator $N$.
\end{lemma}

\begin{proof}
For simplicity, we assume $e > 1$ (the proof is little more technical
when $e = 1$ and is left to the reader in this case). Assume by
contradiction that such an $N : \calM \to \calM$ exists. Put $x_1 =
N(e_1)$ and $x_2 = N(e_2)$. The relation $N \phi = p \phi N$ implies the
following equalities :
$$(S) : \left\{ \begin{array}{l}
p x_1 + p u x_2 = \phi(x_1) + p u e_2 \\
u x_1 + x_2 = p \phi(x_2) + u e_1. \end{array} \right. $$
For all integer $n$, denote by $J_n$ the topological closure of the ideal 
of $S$ generated by all
$\frac{u^i}{q(i)!}$ for $i \geq n$, where $q(i)$ is the quotient in the
Euclidean division of $i$ by $e$. The first equation of $(S)$ shows
that $\phi(x_1) \in J_1 \calM$, and consequently $x_1 \in J_1 \calM$.
From $\phi(J_1) \subset J_p$, we deduce $\phi(x_1) \in J_p \calM$. By
the same way, it follows from the second equation of $(S)$ that $x_2
\equiv p \phi(x_2) \pmod {J_1}$. Since $S/J_1 \simeq W$, this congruence
proves that $x_2 \in J_1 \calM$ and then, as before, $\phi(x_2) \in J_p
\calM$. Resolving $(S)$, we obtain :
$$x_1 \equiv -\frac{u^2}{1-u^2} e_1 + \frac u{1-u^2} e_2 \pmod {J_p
\calM}$$
which gives $\phi(x_1) \equiv u^p e_2 \pmod{J_{p+1} \calM}$. Hence,
$\phi(x_1)$ is not divisible by $p$ in $S$ (here, we use $e > 1$). But,
on the other hand, the first equation of $(S)$ shows directly that
$\phi(x_1)$ have to be divisible by $p$. This is a contradiction.
\end{proof}

\noindent

Briefly, we have an inclusion $\Modphin r S \subset \Modphi r S$ but it is
always strict if $r > 1$. We call $G_\infty$-representations arising
from objects of $\Modphi r S$ \emph{quasi-semi-stable} representations.
Note that if $V$ is a lattice in a semi-stable representation of $G_K$,
its restriction to $G_\infty$ is quasi-semi-stable\footnote{The converse
is not true in general. In fact, there exists a full subcategory of
$\Modphi r S$, whose objects are called \emph{quasi-strongly divisible
lattices}, which is anti-equivalent to the category of
$G_\infty$-lattices in semi-stable representations. See \cite{liu} for
details.}.

\section[The category $\text{\rm Mod}^{r,\phi,N}_{/\si_\infty}$]{The
category $\Modphi r {\si_\infty}$}

The case of quasi-semi-stable representations is simpler because we lay
out an alternative category (defined by Breuil and studied by Kisin) to
describe them. In this section, we give definitions and basic properties
of this category and we prove that it is equivalent to the category of
Breuil modules.

\subsection{Definitions and basic properties}

We relax the condition $r < p-1$ and assume only $r \in \{0, 1, 2, 3,
\ldots, \infty\}$.

\subsubsection*{Objects of linear algebra}

Put $\si = W[[u]]$ and endow it with a Frobenius $\phi : \si \to \si$
defined by: $$\phi\pa{\sum_{n=0}^\infty a_n u^n} = \sum_{n=0}^\infty
\sigma(a_n) u^{pn}.$$ Put also $\si_1 = \si/p\si = k[[u]]$.
As in \S \ref{subsec:breuilmod}, we define some categories of modules
over $\si$. First, the ``big'' category $\pModphi r \si$: if $r$
is finite, its objects are the $\si$-modules $\frakM$ equipped with a
$\phi$-semi-linear endomorphism $\phi : \frakM \to \frakM$ such that
\begin{equation}
\label{eq:condphi}
E(u)^r \frakM \subset \brac{\im \phi}
\end{equation}
where $\brac{\im \phi}$ denotes the $\si$-submodule of $\frakM$
generated by the image of $\phi$. If $\phi^\star \frakM = \si
\otimes_{(\phi), \si} \frakM$, the previous condition is equivalent to
ask the cokernel of $\id \otimes \phi : \phi^\star \frakM \to \frakM$ to
be killed by $E(u)^r$. If $r = \infty$, we ask condition
(\ref{eq:condphi}) for a non fixed integer $r$: in this way,
$\pModphi \infty \si$ is just the union (in an obvious sense) of all
categories $\pModphi r \si$ for $r$ finite. Morphisms in $\pModphi r \si$
are just $\si$-linear morphisms that commute with Frobenius.

Now, we define full subcategories of $\pModphi r \si$. The category
$\Modphi r \si$ (resp. $\Modphi r {\si_1}$) gathers all objects $\frakM \in
\pModphi r \si$ free of finite rank over $\si$ (resp. over $\si_1$),
whereas $\Modphi r {\si_\infty}$ is the smallest subcategory of
$\pModphi r \si$ containing $\Modphi r {\si_1}$ and stable under
extensions\footnote{An sequence of objects of $\pModphi r \si$ is said
exact if it is exact as a sequence of $\si$-modules.}. For simplicity,
we also define the category $\pModphi r {\si_\infty}$ as the full
subcategory of $\pModphi r \si$ gathering all objects killed by a power of
$p$. Obviously we have $\Modphi r {\si_\infty} \subset
\pModphi r {\si_\infty}$. The following proposition summarizes basic
properties of these modules.

\begin{prop}
\label{prop:basickisin}
\begin{itemize}
\item[(i)] Let $\frakM \in \Modphi r {\si_\infty}$. Then $\id \otimes \phi
: \phi^\star \frakM \to \frakM$ is injective.

\item[(ii)] Let $\frakM$ be an object of $\pModphi r \si$. Then $\frakM$
is in $\Modphi r {\si_\infty}$ if and only if it is of finite type over
$\si$, it have no $u$-torsion and it is killed by a power of $p$.

\item[(iii)] The category $\Modphi r {\si_\infty}$ is stable under kernels
and images.
\end{itemize}
\end{prop}

\begin{proof}
See \S 2.3 of \cite{liu2}.
\end{proof}

\medskip

Furthermore, there is a functor $M_{\si_\infty} : \Modphi r {\si_\infty}
\to \pModphi r S$. It is defined as follows. Let $\frakM$ be an object
of
$\Modphi r {\si_\infty}$. As an $S$-module, $M_\si (\frakM) = S
\otimes_{(\phi), \si} \frakM$ where the subscript ``$(\phi)$'' means
that $S$ is considered as a $\si$-module \emph{via} the composite $\si
\to S \to S$, the first map being the canonimal map and the second the
Frobenius $\phi$. The Frobenius on $\frakM$ induces a $S$-linear map
$\id \otimes \phi : \calM \to S \otimes_\si \frakM$. We then define
$\Fil^r \calM$ by the formula
$$\Fil^r \calM = \acco{x \in \calM, \, (\id \otimes \phi)(x) \in \Fil^r
S \otimes_\si \frakM \subset S \otimes_\si \frakM)}.$$
The map $\phi_r$ is given by the following composite:
$$\xymatrix {
\Fil^r \calM \ar[r]^-{\id \otimes \phi} & \Fil^r S \otimes_\si \frakM
\ar[r]^-{\phi_r \otimes \id} & \calM.}$$
Identical constructions give rise to an other functor $M_\si : \Modphi
r \si \to \pModphi r S$.

\begin{prop}
The functor $M_{\si_\infty}$ (resp. $M_\si$) takes values in $\Modphi r
{S_\infty}$ (resp. $\Modphi r S$). Moreover, both functors are exact and
fully faithful.
\end{prop}

\begin{proof}
The case $r=1$ is done in proposition 1.1.11 of \cite{kisin04}. The same
proof works for any $r$.
\end{proof}

\begin{prop}
\label{prop:quotmodsi}
Let $\frakM' \subset \frakM$ be two objects of $\Modphi r \si$ such that
$\frakM' \otimes_{\Z_p} \Q_p \simeq \frakM \otimes_{\Z_p} \Q_p$. Then,
the quotient $\frakM'' = \frakM / \frakM'$ is an object of $\Modphi r
{\si_\infty}$. Moreover, the sequence
$$0 \to M_\si (\frakM') \to M_\si (\frakM) \to M_{\si_\infty} (\frakM'')
\to 0$$
is exact.
\end{prop}

\begin{proof}
The first point is proved in proposition 2.3.2 of \cite{liu2}.
For the second point, the proof is the same as for the exactness of
$M_{\si_\infty}$.
\end{proof}

\subsubsection*{Functors to Galois representations}

We recall the construction of the functor $\pT_{\si_\infty}$ from
$\pModphi r {\si_\infty}$ to the category of $\Z_p$-representations of
$G_\infty$.
First, we define several rings. Put $R = \varprojlim \O_{\bar K}/p$
where the transition maps are given by Frobenius. There is a unique
surjective map $\theta : W(R) \to \widehat{\O_{\bar K}}$ to the $p$-adic
completion $\widehat{\O_{\bar K}}$ of $\O_{\bar K}$, which lifts the
projection $R \to \O_{\bar K}/p$ onto the first factor. Recall that
we have fixed a sequence $\pi_n)_{n \geq 0}$ of compatible $p^n$-th
root of $\pi$. It defines an element of $R$ and we
denote by $[\underline \pi]$ its Teichm\"{u}ller representative. We have an
embedding $\si \to W(R)$, $u \mapsto [\underline \pi]$ which is
compatible with Frobenius.

Let $\O_\E$ be the $p$-adic completion of $\si[1/u]$. It is a discrete
valuation ring with residue field $k((u))$. Put $\E = \Frac \O_\E$. The
embedding $\si \to W(R)$ extends to an embedding $\E \to W(\Frac R)$.
Let $\E^\ur$ the maximal unramified extension of $\E$ included in
$W(\Frac R)[1/p]$ and $\O_{\E^\ur}$ its ring of integers. Since $W(\Frac
R)$ is algebraically closed (see \cite{fontaine-fest}, \S A.3.1.6), the
residue field $\O_{\E^\ur} / p$ is isomorphic to $k((u))^\sep$, a
separable closure of $k((u))$. We will consider the tensor product
$\O_{\E^\ur} \otimes_{\Z_p} \Q_p/\Z_p = \E^\ur / \O_{\E^\ur}$. It is an
object of $\pModphi r \si$ endowed with an action of $G_\infty$.

Finally, the functor $\pT_{\si_\infty}$ is defined by the formula
$$\pT_{\si_\infty} (\frakM) = \hom_{\pModphi r \si} (\frakM, \O_{\E^\ur}
\otimes_{\Z_p} \Q_p/\Z_p)$$
for each $\frakM \in \pModphi r {\si_\infty}$. We call $T_{\si_\infty}$
the restriction of $\pT_{\si_\infty}$ to the subcategory $\Modphi r
{\si_\infty}$. If $\frakM \in \Modphi r {\si_1}$, the expression of
$T_{\si_\infty} (\frakM)$ can be simplified as follows:
$$T_{\si_\infty} (\frakM) = \hom_{\pModphi r \si} (\frakM, \O_{\E^\ur} / p) =
\hom_{\pModphi r \si} (\frakM, k((u))^\sep).$$

\begin{prop}
\label{prop:tsiinf}
The composite $T_\qcris \circ M_{\si_\infty}$ is $T_{\si_\infty}$ and it
is an exact functor.

If $\frakM \in \Modphi r {\si_1}$ is free of rank $d$ over
$\si_1$, then $T_{\si_\infty} (\frakM)$ is a vector space of dimension
$d$ over $\F_p$.
\end{prop}

\begin{proof}
It has been proved in \S B.1.8.4 and \S A.1.2 in \cite{fontaine-fest}.
\end{proof}

\begin{lemma}
\label{lem:reduit}
Let $\frakM \in \Modphi r {\si_\infty}$. Then
$\bigcap_{f \in T_{\si_\infty}(\frakM)} \ker f = 0$.
\end{lemma}

\begin{proof}
First, we show the lemma for $\frakM \in \Modphi r {\si_1}$. Put $\frakK =
\bigcap_{f \in T_{\si_\infty}(\frakM)} \ker f$. Since $u$ in invertible
in $k((u))^\sep$, the quotient $\frakM / \frakK$ have no $u$-torsion and
by proposition \ref{prop:basickisin} (ii), it is an object of $\Modphi r
{\si_1}$.
Furthermore, by definition of $\frakK$, the map $\frakM \to
\frakM/\frakK$ induces a bijection $T_{\si_\infty} (\frakM/\frakK) \to
T_{\si_\infty} (\frakM)$. By proposition \ref{prop:tsiinf}, modules
$\frakM/\frakK$ and $\frakM$ have same rank and hence $\frakK = 0$ as
required.

It remains to prove that if $0 \to \frakM' \to \frakM \to \frakM'' \to
0$ is an exact sequence in $\Modphi r {\si_\infty}$ and if the conclusion
is correct for $\frakM'$ and $\frakM''$, then it is also correct for
$\frakM$. Let $x \in \frakM$ such that $f(x) = 0$ for all $f \in
T_{\si_\infty} (\frakM)$. If $y \in \frakM''$ is the image of $x$, we
have $g(y) = 0$ for all $g \in T_{\si_\infty} (\frakM)$. Thus by
assumption $y = 0$. Hence $x \in \frakM'$. Let $g \in T_{\si_\infty}
(\frakM')$. By exactness of $T_{\si_\infty}$ (proposition
\ref{prop:tsiinf}), $g$ can be extended to a map $f \in T_{\si_\infty}
(\frakM)$. Using the assumption, we get $g(x) = 0$ and finally $x=0$.
\end{proof}

\begin{cor}
\label{cor:sifaithful}
The functor $T_{\si_\infty}$ is faithful.
\end{cor}

\subsection{An equivalence of categories}

The aim of this subsection is to prove the following theorem.

\begin{theo}
\label{theo:equivfree}
Assume $r < p-1$.
The functor $M_\si : \Modphi r \si \to \Modphi r S$ is an equivalence of
categories.
\end{theo}

The full faithfulness was already seen.
Hence it remains to prove the essential surjectivity. Let $\calM \in
\Modphi r S$ and denote by $d$ its rank over $S$. The heart of the proof
is the following technical lemma.

\begin{lemma}
With previous notations, there exists $\alpha_1, \dots, \alpha_d \in
\Fil^r \calM$ and a basis $e_1, \dots , e_d$ of $\calM$ such that $e_i =
\frac 1{c^r} \phi_r(\alpha_i)$, $(\alpha_1, \ldots, \alpha_d) = (e_1 ,
\ldots , e_d) B$ with $B$ a $d \times d$ matrix with coefficients in
$\si$ and
\begin{equation}
\label{eq:genfil}
\Fil^r \calM = \sum_{i=1}^d S\alpha _i + \Fil^pS \calM.
\end{equation}
\end{lemma}

\begin{proof}
If $R$ is a ring, we denote by $M_d(R)$ the algebra of $d \times d$
matrices with coefficients in $R$.

We first show that we can inductively construct $(\alpha^{(n)}_1, \dots,
\alpha_d^{(n)}) \in \Fil^r \calM$ such that
\begin{enumerate}
\item $(e^{(n)}_1 , \dots , e^{(n)}_d)=
c^{-r} \phi_r(\alpha^{(n)}_1 , \dots , \alpha_d^{(n)})$ is a
basis of $\calM$;
\item
there exists matrices  $B^{(n)} \in M_d(\si)$ and $C^{(n)} \in M_d(p^{n}
\Fil^{n + p} S)$ such that $(\alpha^{(n)}_1, \dots, \alpha_d^{(n)})=
(e^{(n)}_1 , \dots , e_d^{(n)})(B^{(n)} + C^{(n)})$.
\end{enumerate}
For $n=0$, the result is a consequence of the (easy part of the) lemma
4.1.1 of \cite{liu}. Note also that property (\ref{eq:genfil}) is satisfied
with $\alpha_i^{(0)}$ instead of $\alpha_i$.
Now, assume that the $\alpha^{(n)}_i$'s are build. We put
\begin{equation}\label{eq:2}
(\alpha^{(n+1)}_1 ,
\dots , \alpha_d^{(n+1)}) =(e^{(n)}_1, \dots , e^{(n)} _d ) B^{(n)}.
\end{equation}
First note that
\begin{eqnarray*}
 (e_1^{(n+1)},\dots, e_d ^{(n+1)})&=& c^{-r} \phi_r(\alpha^{(n+1)}_1 ,
\dots , \alpha_d^{(n+1)})\\
&=&c^{-r} \phi_r((\alpha^{(n)}_1 , \dots , \alpha_d^{(n)})-
(e^{(n)}_1,\dots , e_d ^{(n)}) C^{(n)} ))\\ &=& (e^{(n)}_1 , \dots ,
e_d^{(n)} ) (I -D^{(n)})
\end{eqnarray*}
where $c ^{-r} \phi_r ((e^{(n)}_1,\dots , e_d ^{(n)}) C^{(n)})=
(e^{(n)}_1,\dots , e_d ^{(n)})D^{(n)}$.

Now we claim that $p^{\lambda_n+n}$ divides $D^{(n)}$ where $\lambda_n =
n+p - r -[\frac{n+p}{p-1}]$. Recall that for all $s \in \Fil^rS$ and
$x\in \calM$ we have $\phi_r(sx)= c^{-r}\phi_r(s) \phi_r(E(u)^rx)$.
Moreover, by assumption, $C^{(n)} \in M_d(p^n \Fil^{n+p} S)$. So to
prove the claim it suffices to show that $v_p (\phi_r(s)) \geq
\lambda_n$ for all $s \in \Fil^{n+p}S$. Since $s$ can be always
represented by
$$s = \sum _{m= n+p}^\infty a_m (u)\frac{E(u)^m}{m!}, \ a_m (u) \in
W[u], \ a_m(u) \to 0 \quad p\text{-adically}$$
and $\phi(E(u))= p c$, we reduce the proof to show that
$$ m - v _p (m !)-r > n+p - r- \frac{n+p}{p-1} \ \text{ for any }\ m
\geq n+p$$
which is clear, using $v_p (m!) < \frac{m}{p-1}$.

It is easy to check $\lambda_n \geq 1$. Since $p ^{\lambda_n+n} |
D^{(n)}$, $(I- D^{(n)})$ is invertible and $(e_1^{(n+1)},\dots, e_d
^{(n+1)})$ is a basis of $\calM$. Now by (\ref{eq:2}), we have
$$ (\alpha^{(n+1)}_1 , \dots , \alpha_d^{(n+1)})=(e^{(n)}_1, \dots ,
e^{(n)} _d ) B^{(n)}= (e_1^{(n+1)},\dots, e_d^{(n+1)})
(I-D^{(n)})^{-1}B^{(n)}.$$
Put $A = (I-D^{(n)})^{-1}B^{(n)}$. To achieve the induction, it
remains to write $A = B^{(n+1)} + C^{(n+1)}$ with $B^{(n+1)} \in
M_d(\si)$ and $C^{(n+1)} \in M_d(p^{n+1} \Fil^{n+1+p} S)$. For that,
write $D^{(n)} = p^{\lambda_n+n} E^{(n)}$ and
$$E^{(n)} = \sum_{i=0}^ {n+p} b_i(u) \frac{E(u)^i}{i !} + \sum_{i =
n+p+1}^\infty b_i(u) \frac{E(u)^i}{i !} = E^{(n)}_1 + E^{(n)}_2$$
with $b_i(u) \in W[u]$. A simple computation on valutation gives
$p^{\lambda_n+n} i ! \in \Z_p$ for all $i \leq n+p$.
Thus $D_1^{(n)} = p^{\lambda_n + n} E_1^{(n)} \in M_d(\si)$. The
conclusion then follows by expanding the series
$$A = \sum_{i=0}^\infty (D_1^{(n)} + D_2^{(n)})^i B^{(n)}$$
where $D_2^{(n)} = p^{\lambda_n+n} E_2^{(n)} \in M_d(p^{n+1}
\Fil^{n+1+p} S)$.

\medskip

To complete the proof of the lemma, remark that equation (\ref{eq:2})
implies
\begin{equation}
\label{eq:prop3}
(\alpha^{(n+1)}_1 ,\dots , \alpha_d^{(n+1)}) - (\alpha_1^{(n)}
, \dots , \alpha_d^{(n)}) = - (e^{(n)}_1, \dots , e^{(n)} _d) C^{(n)}
\end{equation}
and hence the convergence of all $\alpha_i^{(n)}$ because $p^n$ divides
$C^{(n)}$. The convergence of all $e_i^{(n)}$ and then those of matrices
$B^{(n)}$ follows. If $\alpha_i$ (resp. $B$) is the limit of
$\alpha_i^{(n)}$ (resp. $B^{(n)}$), we have $\phi_r (\alpha_1, \ldots,
\alpha_d) = c^{-r} (e_1, \ldots, e_d)$ and $(\alpha_1, \ldots, \alpha_d)
= (e_1, \ldots, e_d)B$ with $B \in M_d(\si)$.

It remains to check property (\ref{eq:genfil}). For that, we can reduce
modulo $p$ and then, the conclusion follows from the congruences
$\alpha_i \equiv \alpha_i^{(0)} \pmod p$.
\end{proof}

Now, it is quite easy to achieve the proof of theorem
\ref{theo:equivfree}. First, we show that there exists $A \in M_d(\si)$
such that $BA = E(u)^r I$. Indeed, since $E(u)^r e_i \in \Fil^r \calM$
for all $i$, the condition (\ref{eq:genfil}) implies that there exists
matrices $A'$, $C'$ such that $BA' +C' = E(u)^rI$ and $C' \in M_d(\Fil^p
S)$. Writing $A' = A'_0 + A'_1$ with $A'_0 \in M_d(W[u])$ and $A'_1 \in
M_d(\Fil^pS)$, we may assume $A' \in M_d(W[u])$. Then $C' = E(u)^r I - BA'$
has coefficients in $\si \cap \Fil^p S$. Therefore, $C'= E(u)^p C$ with
$C \in M_d(\si)$. Now $BA'= E(u)^r (I - E(u)^{p-r}C)$ and $A = A'(I -
E(u)^{p-r}C)^{-1} \in M_d(\si)$ is appropriate.

Finally, it is easy to check that $\frakM = \si f_1 \oplus \cdots \oplus
\si f_d$ endowed with $\phi$ defined by $\phi(f_1, \ldots, f_d) =
(f_1, \ldots, f_d)A$ is a preimage of $\calM$ under $M_\si$. This proves
the theorem.

\subsection{Consequences}

The first consequence is the extension of the equivalence on torsion
objects.

\begin{theo}
\label{theo:equivtors}
Assume $r < p-1$.
The functor $M_{\si_\infty} : \Modphi r {\si_\infty} \to \Modphi r {S_\infty}$
is an equivalence of categories.
\end{theo}

\begin{proof}
It remains to show the essential surjectivity.
Let $\calM$ be an object of $\Modphi r {S_\infty}$. By theorem V.2.a
of \cite{caruso-thesis}, there exists two objects $\hat \calM$ and 
$\hat \calM'$ in $\Modphi r {\si_\infty}$, together with an exact sequence 
$0 \to \hat \calM' \to \hat \calM \to \calM \to 0$ in $\pModphi r S$. Now, 
by theorem
\ref{theo:equivfree}, we can find $\hat \frakM$ and $\hat \frakM'$ two
objects of $\Modphi r \si$ such that $M_\si(\hat \frakM) = \hat \calM$ and
$M_\si (\hat \frakM') = \hat \calM'$. We can also find a map $f : \hat
\frakM' \to \hat \frakM$ inducing the canonical inclusion $\hat
\calM' \to \hat \calM$. The map $F = T_\si(f)$ is an injective
application between two free $\Z_p$-modules of same (finite) rank. Consequently,
there exists $G : T_\si(\hat \frakM') \to T_\si(\hat \frakM)$ such that
$F \circ G = G \circ F = p^n \id$ for an integer $n$. By full
faithfulness of $T_\si$, there exists a map $g : \hat \frakM \to \hat
\frakM'$ satisfying $f \circ g = g \circ f = p^n \id$. It follows that
$f \otimes_{\Z_p} \Q_p$ is bijective. Then, we can apply proposition
\ref{prop:quotmodsi}: $\frakM = \hat \frakM / \hat \frakM'$ is in
$\Modphi r {\si_\infty}$ and $M_{\si_\infty} (\frakM) = \calM$. The
theorem follows.
\end{proof}

\begin{prop}
Assume $r < p-1$ and choose $M_{S_\infty}$ a quasi-inverse of
$M_{\si_\infty}$. If $f : \calM
\to \calM'$ is an injective (resp. surjective) morphism in $\Modphi r
{S_\infty}$, then $M_{S_\infty} (f)$ is also. Moreover, the functor
$M_{S_\infty}$ is exact.
\end{prop}

\begin{proof}
Let $f : \calM \to \calM'$ be a morphism in
$\Modphi r {S_\infty}$. Put $\frakM = M_{S_\infty} (\calM)$, $\frakM' =
M_{S_\infty} (\calM')$ and $g = M_{S_\infty} (f)$.

Assume $f$ injective and denote by $\frakK$ the kernel of $g$. By
proposition \ref{prop:basickisin} (iii), we have $\frakK \in \Modphi r
{\si_\infty}$. Put $\calK = M_{\si_\infty} (\frakK)$. Let $h : \calK \to
\calM$ the image under $M_{\si_\infty}$ of the inclusion $\frakK \to
\frakM$. The composite $f \circ h$ is zero and since $f$ is injective,
$h = 0$. By faithfulness, the inclusion $\frakK \to \frakM$ vanishes,
and consequently $\frakK = 0$ and $g$ is injective.

Now suppose $f$ surjective and denote by $\frakC$ the cokernel of $g$.
Then $S \otimes_{(\phi), \si} \frakC = 0$. By reducing modulo $p$, we
get $S_1 \otimes_{(\phi), \si_1} \frakC/p\frakC = 0$. Since
$\frakC/p\frakC$ is a module of finite type over the principal ring
$k[[u]]$, it is a direct sum of some $k[[u]]$ or $k[[u]]/u^n$ for a
suitable integers $n$. By computing the tensor product, it follows that
the only solution is $\frakC / p \frakC = 0$, \emph{i.e} $\frakC =
p\frakC$. Since $\frakC$ is finitely generated, Nakayama's lemma gives
$\frakC = 0$ as required.

For the exactness, take $0 \to \calM' \to \calM \to \calM'' \to 0$ an
exact sequence in $\Modphi r {S_\infty}$. We know that $M_{S_\infty}
(\calM) \to M_{S_\infty} (\calM'')$ is surjective. Call $\frakK$ its
kernel: it is an object of $\Modphi r {\si_\infty}$ and we have an exact
sequence $0 \to \frakK \to M_{S_\infty} (\calM) \to M_{S_\infty}
(\calM'') \to 0$. Applying the exact functor $M_{\si_\infty}$, we see
that $M_{\si_\infty} (\frakK)$ is the kernel of $\calM \to \calM''$.
Hence, it is isomorphic to $\calM'$ and we are done.
\end{proof}

\noindent
{\it Remark.} Although the functor $M_{\si_\infty}$ is exact, the
implication ($f$ injective) $\Rightarrow$ ($M_{\si_\infty} (f)$
injective) is not true if $er \geq p-1$. Here is a counter-example. Take
$\frakM = \si_1$ with $\phi(1) = 1$, $\frakM' = \si_1$ with $\phi(1) =
u^{p-1}$ and $f : \frakM' \to \frakM$, $1 \mapsto u$. It is injective.
However, $\calM = M_{\si_\infty}$ is just $S_1$ endowed with $\Fil^r
S_1$ and the canonical $\phi_r$. On the other hand, $\calM' = S_1$,
$\Fil^r \calM' = u^{er-p+1} \calM'$ and $\phi_r(u^{er-p+1}) = (-1)^r$.
The map $M_{\si_\infty} (f)$ is the multiplication by $u^p$ and sends
$u^{(e-1)p}$ to $0$; hence it is not injective.

\medskip

\begin{cor}
Assume $r < p-1$.
Functors $T_\qcris$ on $\Modphi r {S_\infty}$ and $T_\st$ on $\ModphiN r
{S_\infty}$ are faithful.
\end{cor}

\begin{proof}
For $T_\qcris$, it is a direct consequence of corollary
\ref{cor:sifaithful} and theorem \ref{theo:equivtors}.

Let $f : \calM \to \calM'$ be a morphism in $\ModphiN r {S_\infty}$. It
can be seen as a morphism in $\Modphi r {S_\infty}$ and we have
$T_\qcris(f) = T_\st(f)$. If this morphism vanishes, then $f$ have also
to vanish thanks to the faithfulness of $T_\qcris$. This proves the
corollary.
\end{proof}

\begin{theo}
Assume $r < p-1$.
Let $\calM' \subset \calM$ be two objects of $\Modphi r S$ such that
$\calM' \otimes_{\Z_p} \Q_p \simeq \calM \otimes_{\Z_p} \Q_p$ and
$\Fil^r \calM' \otimes_{\Z_p} \Q_p \simeq \Fil^r \calM \otimes_{\Z_p}
\Q_p$. Then the quotient $\calM / \calM'$ is an object of $\Modphi r
{S_\infty}$. Furthermore every object of $\Modphi r {S_\infty}$ can be
written in this way.
\end{theo}

\begin{proof}
For the first part of the theorem, we use a similar argument as in the
proof of theorem \ref{theo:equivtors}. Let $\frakM' \to \frakM$ an
antecedent of the inclusion $\calM' \to \calM$. We first show that
$\frakM' \otimes_{\Z_p} \Q_p \simeq \frakM \otimes_{\Z_p} \Q_p$, and
then by using proposition \ref{prop:quotmodsi}, we get $M_{\si_\infty}
(\frakM / \frakM') = \calM / \calM'$.

The second part is again theorem V.2.a of \cite{caruso-thesis}.
\end{proof}

\noindent
{\it Remark.} The condition $\Fil^r \calM' \otimes_{\Z_p} \Q_p \simeq
\Fil^r \calM \otimes_{\Z_p} \Q_p$ is equivalent to $\Fil^r \calM' =
\calM' \cap \Fil^r \calM$. Indeed, if $x \in \calM' \cap \Fil^r \calM$
then $x \in \Fil^r \calM' \otimes_{\Z_p} \Q_p = \Fil^r \calM
\otimes_{\Z_p} \Q_p$ and $p^n x \in \Fil^r \calM'$ for a certain integer
$n$. Since, by definition, $\calM' / \Fil^r \calM'$ have no $p$-torsion,
we must have $x \in \Fil^r \calM'$. The controverse is easy.

\subsection{Duality}
\label{subsec:duality}

In \cite{liu2}, \S 3.1, one of the author has defined a duality on
$\Modphi r {\si_\infty}$ for all $r < \infty$. It consists in an exact
functor $\frakM \mapsto \frakM^\vee$. Let's recall its definition and
properties. For $\frakM \in \Modphi r {\si_\infty}$, we put $\frakM^\vee
= \hom_\si (\frakM, \si \otimes_{\Z_p} \Q_p/\Z_p)$. We then have a
natural pairing :
$$\brac{\cdot, \cdot} : \frakM \times \frakM^\vee \to \si \otimes_{\Z_p}
\Q_p/\Z_p.$$
The Frobenius $\phi^\vee$ on $\frakM^\vee$ is defined by the equality
$$\brac{\phi(x), \phi^\vee(y)} = c_0^{-r} E(u)^r \phi(\brac{x,y})$$
(for all $x \in \frakM$ and $y \in \frakM^\vee$) where $c_0 =
\frac{E(0)} p \in W^\star$ and the latest $\phi$ is gotten from the
usual operator on $\si$.

Here are main properties of the duality. We have a natural isomorphism
$(\frakM^\vee)^\vee \simeq \frakM$, and a compatibility between duality
and $T_{\si_\infty}$ given by the following functorial isomorphism:
\begin{equation}
\label{eq:tsidual}
T_{\si_\infty} (\frakM^\vee) \simeq T_{\si_\infty} (\frakM)^\vee (r).
\end{equation}
where ``$(r)$'' is for the Tate twist.

\bigskip

In \cite{caruso-thesis}, chapter V, one of the author (not the same) has
defined a duality on $\Modphi r {S_\infty}$ for $r < p-1$. If $\calM$ is
an object of this category, we put $\calM^\vee = \hom_S (\calM, S
\otimes_{\Z_p} \Q_p/\Z_p)$, $\Fil^r \calM^\vee = \{ f \in \calM^\vee, \,
f(\Fil^r \calM) \subset \Fil^r S \otimes_{\Z_p} \Q_p/\Z_p\}$ and if $f
\in \Fil^r \calM^\vee$,
$\phi_r^\vee(f)$ is defined as the unique map making commutative the
following diagram:
$$\xymatrix {
\Fil^r \calM \ar[r]^-{\phi_r} \ar[d]_-{f} & \calM \ar[d]^-{\phi_r^\vee
(f)} \\
\Fil^r S \otimes_{\Z_p} \Q_p/\Z_p \ar[r]^-{\phi_r} & S \otimes_{\Z_p}
\Q_p/\Z_p }$$
Now, consider $\frakM \in \Modphi r {\si_\infty}$ (always with $r <
p-1$). Put:
$$\lambda = \prod_{n=1}^\infty \phi^n\pa{\frac{E(u)}{pc_0}} \in S.$$
and define the following canonical isomorphism:
$$M_{\si_\infty} (\frakM^\vee) \to M_{\si_\infty}(\frakM)^\vee, \quad
s \otimes f \mapsto \frac 1 {\lambda^r} \: s f.$$
A direct calculation gives $\phi(\lambda) = \frac c{\phi(c_0)} \lambda$,
which implies that the previous isomorphism is compatible with $\phi$,
and hence a morphism in $\Modphi r {S_\infty}$. We deduce the following:

\begin{cor}
\label{cor:duality}
Assume $r < p-1$. For any $\calM \in \Modphi r {S_\infty}$, there exists
a natural isomorphism $\calM \to (\calM^\vee)^\vee$ and a natural
isomorphism:
$$T_{\qcris} (\calM^\vee) \simeq T_{\qcris} (\calM)^\vee (r).$$
\end{cor}

\noindent
{\it Remarks.} Corollary \ref{cor:duality} is proved (with different
methods) in \cite{caruso-thesis} under the assumption $er < p-1$ or
$r=1$.

In \emph{loc. cit.}, definition of duality is extended to the category
$\ModphiN r {S_\infty}$: the operator $N^\vee$ on $\calM^\vee$ is
defined by the formula $N^\vee(f) = N \circ f - f \circ N$ (where $N$ is
the given operator on $\calM$). Using isomorphism (\ref{eq:tsttqst}), we
directly obtain a version of corollary \ref{cor:duality} in this new
situation.

\section[A construction on $\text{\rm Mod}^{r,\phi,N}_{/\si_\infty}$]{A
construction on $\Modphi r {\si_\infty}$}
\label{sec:const}

This section is devoted to give a proof of theorem \ref{theo:overS}. We
will use the equivalence stated in theorem \ref{theo:equivtors} to make
constructions with more pleasant modules.

\subsection[The category $\text{\rm 'Mod}^{\phi,N}_{/\O_\E}$]{The
category $\pModphi {} {\O_\E}$}

Let's recall classical results about the classification of
$\Z_p$-representations of $G_\infty$. Denote by $\pModphi {} {\O_\E}$
the category of torsion \'etale $\phi$-modules over $\O_\E$. By
definition, an object of $\pModphi {} {\O_\E}$ is an $\O_\E$-module
$M$ killed by a power of $p$ and equipped with a Frobenius $\phi :
M \to M$ that induces a bijection $\id \otimes \phi :
\phi^\star M \to M$ (where $\phi^\star M = \O_\E
\otimes_{(\phi), \O_\E} \frakM$).

\bigskip

\noindent
{\it Remark.} Since we are only interested in $p$-torsion modules, the
definition does not change if we substitute the ring $\si[1/u]$ to
$\O_\E$ (in other words, we do not need to complete $p$-adically). In
the sequel, we will just work with $\si[1/u]$.

\bigskip

We have a functor $\pT_{\O_\E} : \pModphi {} {\O_\E} \to \Rep_{\Z_p}
(G_\infty)$ defined by
$$\pT_{\O_\E} (M) = \hom_{\pModphi {} {\O_\E}} (M, \O_{\E^\ur}
\otimes_{\Z_p} \Q_p/\Z_p).$$

\begin{theo}
The functor $\pT_{\O_\E}$ is exact and fully faithful.
\end{theo}

\begin{proof}
See \S A.1.2 of \cite{fontaine-fest}.
\end{proof}

Furthermore $\pT_{\si_\infty}$ factors through $\pT_{\O_\E}$ as follows:
if $\pM_{\O_\E} : \Modphi r {\si_\infty} \to \pModphi {} {\O_\E}$ is
defined by $\pM_{\O_\E} (\frakM) = \frakM \otimes_\si \O_\E = \frakM
\otimes_\si \si[1/u]$ (since $E(u)$ is invertible in $\O_\E$, the map
$\id \otimes \phi : \phi^\star [\pM_{\O_\E} (\frakM)] \to \pM_{\O_\E}
(\frakM)$ is bijective), the equality $\pT_{\si_\infty} = \pT_{\O_\E}
\circ \pM_{\O_\E}$ holds. In a slightly different situation,
$\pM_{\O_\E}$ is the functor $j^\star$ of \cite{fontaine-fest}. From
now on, we will use the notation $\frakM[1/u]$ for $\pM_{\O_\E} (\frakM)$.
In \cite{fontaine-fest}, Fontaine defines an adjoint $j_\star$ to his
functor $j^\star$. In the sequel, we will adapt his construction to our
settings.

\subsection{The ordered set $F^r_\si(M)$}

In this subsection, we fix $M \in \pModphi {} {\O_\E}$. Our aim is to
study the structure of the ``set'' of previous images of $M$ under
$\pM_{\O_\E}$. We begin by the following definition:

\begin{deftn}
\label{def:frsi}
Let $\calF^r_\si (M)$ the category whose objects are couples $(\frakM,
f)$ where $\frakM$ is an object of $\Modphi r {\si_\infty}$ and
$f : \frakM[1/u] \to M$ is an isomorphism. Morphisms in $\calF^r_\si
(M)$ are morphisms in $\Modphi r {\si_\infty}$ that are compatible with
$f$.

Let $F^r_\si (M)$ be the (partially) ordered set (by inclusion) of
$\frakM \in \Modphi r {\si_\infty}$ contained in $M$ such that
$\frakM[1/u] = M$.
\end{deftn}

\medskip

The following lemma is easy:

\begin{lemma}
The category $\calF^r_\si(M)$ is equivalent to (the category associated
to) the ordered set $F^r_\si(M)$.
\end{lemma}

\subsubsection*{Supremum and infimum}

\begin{prop}
\label{prop:supinf}
The ordered set $F^r_\si (M)$ has finite supremum and finite infimum.
\end{prop}

\begin{proof}
Obviously, it suffices to prove that for any $\frakM'$ and
$\frakM''$ in $F^r_\si(M)$, $\sup(\frakM', \frakM'')$ and $\inf(\frakM',
\frakM'')$ exist.

For the supremum, it is enough to show that $\frakM = \frakM' +
\frakM''$ (where the sum is computed in $M$) is an object of $\Modphi r
{\si_\infty}$ (it is obvious that $\frakM[1/u] = M$). For this, remark
that since $\frakM'$ and $\frakM''$ satisfy condition
(\ref{eq:condphi}) (defined page \pageref{eq:condphi}), $\frakM$ also.
The conclusion then follows from proposition \ref{prop:basickisin} (ii).

In the same way, for the infimum, we want to prove that $\frakM =
\frakM' \cap \frakM''$ satisfies $\frakM[1/u] = M$ and is in $\Modphi r
{\si_\infty}$. Since $\frakM'$ is finitely generated, there exists an
integer $s$ such that $u^s \frakM' \subset \frakM''$ and the first point
is clear. Now, Let $x \in \frakM$. Because $\frakM'$ and $\frakM''$ are
in $\Modphi r {\si_\infty}$, there exists $x' \in \phi^\star \frakM'$ and
$x'' \in \phi^\star \frakM''$ such that $E(u)^r x = \id \otimes \phi(x')
= \id \otimes \phi(x'')$ (if $r = \infty$, it must be replaced by a
sufficiently large integer). But, by definition, $\id \otimes \phi$
is injective on $\phi^\star M$. It follows that $x' = x'' \in \phi^\star
\frakM$. Consequently, condition (\ref{eq:condphi}) holds for $\frakM$.
Moreover, since $\si$ in noetherian, $\frakM \subset \frakM'$ is
finitely generated over $\si$. Finally, it is obviously killed by a
power of $p$, and without $u$-torsion. Proposition \ref{prop:basickisin}
ends the proof.
\end{proof}

\subsubsection*{Some finiteness property}

\begin{lemma}
\label{lem:boundmax}
Fix $\frakM \in F^r_\si (M)$. There exists an integer $\ell$ (depending
only on $\frakM$) such that $\lg_\si(\frakM' / \frakM) \leq \ell$ for
any $\frakM' \in F^r_\si(M)$ with $\frakM \subset \frakM'$.
\end{lemma}

\begin{proof}
First, we prove by \emph{d\'evissage} that it is sufficient to consider
the case where $M$ is killed by $p$. Denote by $\frakM(p)$ (resp.
$\frakM'(p)$) the kernel of the multiplication by $p$ on $\frakM$ (resp.
$\frakM'$). We have the following commutative diagram:
$$\xymatrix @R=10pt {
0 \ar[r] & \frakM(p) \ar[r] \ar[d] & \frakM \ar[r] \ar[d] &
\frakM/\frakM(p) \ar[r] \ar[d] & 0 \\
0 \ar[r] & \frakM'(p) \ar[r] & \frakM' \ar[r] & \frakM'/\frakM'(p)
\ar[r] & 0 }$$
where both horizontal sequences are exact, and all vertical arrows are
injective. Snake lemma then shows that the sequence $0 \to
\frac{\frakM'(p)}{\frakM(p)} \to \frac{\frakM'}{\frakM} \to
\frac{\frakM'/\frakM(p)}{\frakM'/\frakM(p)} \to 0$ remains exact.
The induction follows.

Since $\id \otimes \phi : \phi^\star \frakM \to \frakM$ is injective
(proposition \ref{prop:basickisin} (i)), the map $\frakM / u \frakM \to
\brac{\im \phi} / u \brac {\im \phi}$ induced by $\phi$ is also
injective. By definition, there exists an integer $s$ such that $E(u)^s
\frakM \subset \brac{\im \phi}$. (If $r$ is finite, we can choose $s =
r$.) It follows the implication
\begin{equation}
\label{eq:implication}
(x \not\in u \frakM) \, \Longrightarrow \, (\phi(x) \not\in u^{es+1}
\frakM).
\end{equation}
Furthermore, there exists an integer $n$ such that $u^n \frakM' \subset
\frakM$. Choose $n$ minimal (not necessary positive). Then, we can find
$x \in \frakM'$ such that $u^{n-1} x \not\in \frakM$. Therefore $u^n x
\in \frakM$ but $u^n x \not\in u\frakM$. By applying implication
(\ref{eq:implication}), we get $\phi(u^n x) \not\in u^{es+1} \frakM$,
then $u^n \phi(x) \not\in u^{1+es-(p-1)n} \frakM$. On the other hand,
$u^n \phi(x) \in u^n \frakM' \subset \frakM$. It follows the inequality
$1 + es - (p-1) n \geq 0$ which gives $n \leq t = E(\frac{es+1}{p-1})$
(here $E$ denotes the integer part). From $u^n \frakM' \subset \frakM$,
we get $u^t \frakM' \subset \frakM$ and the conclusion follows (with
$\ell = t \dim_{k((u))} M$).
\end{proof}

\begin{lemma}
Assume $r < \infty$. There exists an integer $\ell$ (depending only on
$M$) such that $\lg_\si(\frakM' / \frakM) \leq \ell$ for any $\frakM$
and $\frakM'$ in $F^r_\si (M)$ with $\frakM \subset \frakM'$.
\end{lemma}

\begin{proof}
Proof of lemma \ref{lem:boundmax} shows that $\ell$ can be chosen
equal to $\lg_{\O_\E} (M) \times E(\frac{er+1}{p-1})$, which depends
only on $M$.
\end{proof}

\begin{cor}
\label{cor:frsi}
The ordered set $F^r_\si (M)$ always has a greatest element.
Furthermore, if $r < \infty$, $F^r_\si$ is finite and has a smallest
element.
\end{cor}

\noindent
{\it Remark.} Proof of lemma \ref{lem:boundmax} gives an upper bound
for the length of any chain in $F^r_\si(M)$, that is :
$$1 + \lg_{\O_\E} (M) \times E\pa{\frac{er+1}{p-1}}.$$
In particular, if $er < p-1$, the set $F^r_\si (M)$ contains at most
one element. This latest assertion will be used several times in the
sequel.

\subsubsection*{Functoriality}

In view of possible generalizations, we would like to rephrase quickly
previous properties in a more categorical and functorial way.

\begin{prop}
The category $\calF_\si(M)$ has finite (direct) sums and finite products.
\end{prop}

\begin{prop}
The category $\calF_\si(M)$ is noetherian in the following sense: if
$$\xymatrix {
\frakM_1 \ar[r]^-{f_1} & \frakM_2 \ar[r]^-{f_2} & \cdots
\ar[r]^-{f_{n-1}} & \frakM_n \ar[r]^-{f_n} & \cdots } $$
is an infinite sequence of morphisms, all $f_n$ are isomorphisms for
$n$ big enough.

If $r$ is finite, the category $\calF_\si(M)$ is artinian in the
following sense: if $\xymatrix @C=15pt {
\frakM_1 & \ar[l]_-{f_1} \frakM_2 & \ar[l]_-{f_2} \cdots }$
is an infinite sequence of morphisms, all $f_n$ are isomorphisms for
$n$ big enough.
\end{prop}

\begin{prop}
\label{prop:supmorp}
Let $\frakM_1, \ldots \frakM_n$ (resp. $\frakM'_1, \ldots \frakM'_n$)
be objects of $\calF_\si(M)$ (resp. $\calF_\si(M')$). Let $f_i :
\frakM_i \to \frakM'_i$ be morphisms in $\Modphi r {\si_\infty}$.
Put $\frakM = \sup(\frakM_1, \ldots, \frakM_n)$ and $\frakM' =
\sup(\frakM'_1, \ldots, \frakM'_n)$. Then, there exists a unique
map $f : \frakM \to \frakM'$ making commutative all diagrams
$$\xymatrix @R=10pt {
\frakM_i \ar[d] \ar[r]^-{f_i} & \frakM'_i \ar[d] \\
\frakM \ar[r]^-{f} & \frakM' \\
}$$
We put $f = \sup (f_1, \ldots, f_n)$.

Furthermore, the association $(f_1, \ldots, f_n) \mapsto \sup (f_1,
\ldots, f_n)$ is functorial in an obvious sense.
\end{prop}

\begin{proof}
Quite clear after the description of $\sup$ given by the proof of
proposition \ref{prop:supinf}.
\end{proof}

\noindent
{\it Remark.} Of course, the analogous statement with $\inf$ is also
true.

\bigskip

\noindent
{\it Important remark.} Since $\pT_{\O_\E}$ is fully faithful, the
functor $\pM_{\O_\E}$ can be replaced by $T_\si$ in definition
\ref{def:frsi}. Hence, it is possible to define supremum and infimum
without reference to the auxiliary category $\pModphi {} {\O_\E}$.

\subsection{Maximal objects}

In this subsection, we give (and prove) some pleasant properties of
objects arising as the greatest element of one set $F_\si(M)$.

\subsubsection*{The functor $\Max^r$}

\begin{deftn}
\label{def:maxr}
Let $\frakM \in \Modphi r {\si_\infty}$. We define $\Max^r(\frakM)$ to
be the greatest element of $F^r_\si(\frakM[1/u])$. It is endowed with an
homomorphism  $\imax^\frakM : \frakM \to \Max^r(\frakM)$ in the category
$\Modphi r {\si_\infty}$.

An object $\frakM$ of $\Modphi r {\si_\infty}$ is said
\emph{maximal} (in $\Modphi r {\si_\infty}$)\footnote{When the value of
$r$ in clear by the context, we will only say \emph{maximal}.} if the
map $\imax^\frakM$ is an isomorphism.
\end{deftn}

\noindent
{\it Remarks.} By \S B.1.5.3 of \cite{fontaine-fest}, a $\phi$-module over
$\si$ \emph{killed by a power of $p$} satisfies condition
(\ref{eq:condphi}) with $r = \infty$, if and only if $\id \otimes \phi :
\phi^\star \frakM[1/u] \to \frakM[1/u]$ is bijective. It follows that
for any $\frakM \in \Modphi \infty {\si_\infty}$, $\Max^\infty (\frakM)
= j_\star (\frakM[1/u])$ where $j_\star$ is the functor defined in \S
B.1.4 of \emph{loc. cit.}

In general, $\Max^r(\frakM)$ and $\Max^{r+1}(\frakM)$ does not coincide.
For instance, take $r$ such that $er \geq p$ and consider $\frakM = \si
e_1 \oplus \si e_2$ with $\phi(e_1) = u e_1 + u^{er} e_2$ and $\phi(e_2)
= u^p e_1$. Then, $\frakM$ is maximal in $\Modphi r {\si_\infty}$ but
\emph{not} in $\Modphi {r+1} {\si_\infty}$ since the submodule of
$\frakM[1/u]$ generated by $e_1$ and $\frac{e_2} u$ is in $F^{r+1}_\si
(\frakM[1/u])$.

\begin{prop}
The previous definition gives rise to a functor $\Max^r : \Modphi r
{\si_\infty} \to \Modphi r {\si_\infty}$.
\end{prop}

\begin{proof}
We have to prove that any map $f : \frakM \to \frakM'$ induces a map
$\Max^r(\frakM) \to \Max^r(\frakM')$. Let $g = f \otimes_\si \si[1/u]$.
By proposition \ref{prop:basickisin} (iii), $g(\Max^r(\frakM))$ is in
$\Modphi r {\si_\infty}$. Hence $g(\Max^r(\frakM)) \subset
\Max^r(\frakM')$ and we are done.
\end{proof}

\noindent
{\it Remark.}
The collection of homomorphisms $(\imax^\frakM)$ defines a natural
transformation between the identity functor and $\Max^r$.

\bigskip

We now show several properties of the functor $\Max^r$.

\begin{prop}
The functor $\Max^r$ is a projection, that is $\Max^r \circ \Max^r =
\Max^r$. Thus, for any $\frakM \in \Modphi r {\si_\infty}$, the object
$\Max^r(\frakM)$ is maximal.
\end{prop}

\begin{proof}
Just remark that $\Max^r(\frakM)[1/u] = \frakM[1/u]$.
\end{proof}

\begin{prop}
\label{prop:maxleftexact}
The functor $\Max^r$ is left exact.
\end{prop}

\begin{proof}
Let $0 \to \frakM' \to \frakM \to \frakM'' \to 0$ an exact sequence in
$\Modphi r {\si_\infty}$. We have the following commutative diagram:
$$\xymatrix {
0 \ar[r] & \frakM' \ar[r] \ar[d]^-{\imax^{\frakM'}} & \frakM \ar[r]
\ar[d]^-{\imax^\frakM} & \frakM'' \ar[d]^-{\imax^{\frakM''}} \ar[r] & 0 \\
0 \ar[r] & \Max^r(\frakM') \ar[r] \ar@{^(->}[d] & \Max^r(\frakM) \ar[r]
\ar@{^(->}[d] & \Max^r(\frakM'') \ar@{^(->}[d] \\
0 \ar[r] & \frakM'[1/u] \ar[r] & \frakM[1/u] \ar[r] & \frakM''[1/u]
\ar[r] & 0 }$$
where the first line is exact by assumption and the last one is also
exact because of the flatness of $\si[1/u]$ over $\si$. We have to show
that the middle line is exact. Injectivity is obvious.

Let's prove the equality $\Max^r(\frakM') = \Max^r(\frakM) \cap
\frakM'[1/u]$. The inclusion $\subset$ is clear. Now, remark that
$\frakM'_\max = \Max^r(\frakM) \cap \frakM'[1/u]$ is a $\si$-submodule
of
$\frakM'[1/u]$ of finite type, which is stable under $\phi$. Moreover,
consider $x \in \frakM'_\max$. Then, there exists $y \in \phi^\star
\Max^r(\frakM)$ and $z \in \phi^\star \frakM'[1/u]$ such that $E(u)^r x
= \id \otimes \phi(y) = \id \otimes \phi(z)$ (if $r = \infty$, it must
be replaced by a sufficiently large integer). Since $\id \otimes
\phi : \phi^\star \frakM[1/u] \to \frakM[1/u]$ is injective, we have $y
= z \in \phi^\star \frakM'_\max$. Hence $\frakM'_\max$ is an object
of $\pModphi r \si$ and the claimed equality is indeed true. This gives
directly the exactness at middle.
\end{proof}

\noindent
{\it Remark.} Unfortunately, $\Max^r$ is not right exact (even on
$\Modphi r
{\si_1}$) if $er \geq p-1$. For instance, consider $\frakM = \si_1 e_1
\oplus \si_1 e_2$ equipped with $\phi$ defined by $\phi(e_1) = e_1$ and
$\phi(e_2) = u e_1 + u^{p-1} e_2$. Denote by $\frakM'$ the submodule of
$\frakM$ generated by $e_1$. We can easily see that $\frakM$ and
$\frakM'$ are both \emph{maximal} objects of $\Modphi r {\si_1}$. However,
$\frakM / \frakM'$ is isomorphic to $\si_1$ with $\phi(1) = u^{p-1}$. It 
is not maximal since $\frac 1 u \si_1$ is finitely generated and stable
under $\phi$.

\begin{prop}
\label{prop:universal}
Let $\frakM \in \Modphi r {\si_\infty}$. The couple $(\Max^r(\frakM),
\imax^\frakM)$ is characterized by the following universal property:
\begin{itemize}
\item the morphism $T_{\si_\infty} (\imax^\frakM)$ is an isomorphism;
\item for each couple $(\frakM', f)$ where $\frakM' \in \Modphi r
{\si_\infty}$ and $f : \frakM \to \frakM'$ becomes an isomorphism under
$T_{\si_\infty}$, there exists a unique map $g : \frakM' \to
\Max^r(\frakM)$ such that $g \circ f = \imax^\frakM$.
\end{itemize}
\end{prop}

\begin{proof}
The first point is clear.
Take $(\frakM',f)$ as in the proposition. Since the quotient $\frakM /
\Max^r(\frakM)$ is killed by a power of $u$, the map $g$ is uniquely
determinated. On the other hand, by full faithfulness of $\pT_{\O_\E}$,
$f$ induces an isomorphism $\tilde f : \frakM[1/u] \to \frakM'[1/u]$.
Denote by $g$ the restriction of $\tilde f^{-1}$ to $\frakM'$. Since
$\frakM'$ is finitely generated over $\si$, $g(\frakM')$ is also and
hence $g(\frakM') \subset \Max^r(\frakM)$ (by definition of $\Max^r$).
In
other words, $g$ induces a map $\frakM' \to \Max^r(\frakM)$ and it is
easy
to check that $g \circ f = \imax^\frakM$.

It remains to prove that the universal property characterizes
$\Max^r(\frakM)$. But if $\frakM'$ satisfies also the universal
property,
we get two maps $\frakM' \to \Max^r (\frakM)$ and $\Max^r(\frakM) \to
\frakM'$ whose composites must be identity.
\end{proof}

\subsubsection*{The category $\Maxphi r {\si_\infty}$}

\begin{deftn}
We put $\Maxphi r {\si_\infty} = \Max^r(\Modphi r {\si_\infty})$. It is a
full subcategory of $\Modphi r {\si_\infty}$.
\end{deftn}

We now show several pleasant properties of this category.

\begin{prop}
\label{prop:adjonction}
The functor $\Max^r : \Modphi r {\si_\infty} \to \Maxphi r {\si_\infty}$ is
a left adjoint to the inclusion functor $\Maxphi r {\si_\infty} \to
\Modphi r {\si_\infty}$.
\end{prop}

\begin{proof}
Let $f : \frakM \to \frakM'$ a morphism in $\Modphi r {\si_\infty}$ and
assume that $\frakM'$ is maximal. We have to prove that there exists a
unique map $\tilde f : \Max^r(\frakM) \to \frakM'$ such that $\tilde f
\circ \imax^\frakM = f$. The unicity is implied by the following
observation: $\frakM'$ have no $u$-torsion, and $\Max^r(\frakM) /
\frakM$ is cancelled by a power of $u$. For the existence, just remark
that $\tilde f = \Max^r(f)$ is appropriate.
\end{proof}

\begin{theo}
\label{theo:abelian}
The category $\Maxphi r {\si_\infty}$ is abelian. More precisely, if $f :
\frakM \to \frakM'$ is a morphism in $\Maxphi r {\si_\infty}$
\begin{itemize}
\item the kernel of $f$ in the usual sense is an object of
$\Maxphi r {\si_\infty}$ and is the kernel of $f$ in the abelian
category $\Maxphi r {\si_\infty}$ ;
\item the cokernel of $f$ in the usual sense, $\coker f$, is an object
of $\pModphi r {\si_\infty}$ and $\Max^r(\frac{\coker f}{u\text{\rm
-torsion}})$ is the cokernel of $f$ in the abelian category $\Maxphi r
{\si_\infty}$; moreover if $f$ is injective, then $\coker f$ have no
$u$-torsion ;
\item the image (resp. coimage) of $f$ in the usual sense is an object
of $\Modphi r {\si_\infty}$ and its image under the functor $\Max^r$ is
the image (resp. coimage) of $f$ in the abelian category $\Maxphi r
{\si_\infty}$.
\end{itemize}
\end{theo}

\begin{proof}
Let $f : \frakM \to \frakM'$ be a morphism in $\Maxphi r {\si_\infty}$. By
proposition \ref{prop:basickisin} (iii), $\frakK = \ker f$ is in object
of $\Modphi r {\si_\infty}$. It remains to prove that it is maximal.
Denote by $\frakM_\max$ the $\si$-submodule of $\frakM[1/u]$ generated
by $\Max^r(\frakK)$ and $\frakM$. It satisfies condition
(\ref{eq:condphi}) (because $\Max^r(\frakK)$ and $\frakM$ satisfy it)
and hence, by proposition \ref{prop:basickisin} (ii), it is an object of
$\Modphi r {\si_\infty}$ included in $\frakM[1/u]$. Since $\frakM$ is
assumed to be maximal, we get $\frakM_\max \subset \frakM$ and then
$\Max^r(\frakK) \subset \frakM$. It
follows $\Max^r(\frakK) \subset \frakM \cap \frakK[1/u] \subset \frakK$
(for the last inclusion, use $\frakK[1/u] = \ker (f \otimes_\si
\si[1/u])$), and $\Max^r(\frakK) = \frakK$.

With proposition \ref{prop:adjonction}, it is easy to prove that
$\Max^r(\frac{\coker f}{u\text{\rm-torsion}})$ is the cokernel of $f$ 
in $\Maxphi r {\si_\infty}$.
The implication ($f$ injective) $\Rightarrow$ ($\coker f \in \Modphi r
{\si_\infty}$) is showed as in proposition \ref{prop:maxleftexact}. It
remains to prove the last statement. We have already seen that the
usual image of $f$, say $\im f$, is an object of $\Modphi r
{\si_\infty}$ (proposition \ref{prop:basickisin} (iii)). Let $g : \im f
\to \frakM'$ the natural inclusion. We have $\coker g = \coker f$. On
the other hand, since $\Max^r(g)$ is an injective morphism between two
maximal objects, its cokernel have no $u$-torsion. Together with 
$g \otimes_\si \si[1/u] = \Max^r(g) \otimes_\si \si[1/u]$, it implies
$\coker \Max^r(g) = \frac{\coker f}{u\text{\rm-torsion}}$. Now, applying 
the left-exact functor $\Max^r$ (see proposition \ref{prop:maxleftexact})
to the exact sequence $0 \to \Max^r(\im f) \to \frakM' \to 
\frac{\coker f}{u\text{\rm-torsion}} \to 0$, we get $\Max^r(\im f) =
\ker(\frakM' \to \frakC)$ where $\frakC = \Max^r(\frac{\coker f}
{u\text{\rm-torsion}}$. Statement about image is then proved.

Finally, by definition, the usual coimage (resp. coimage in $\Maxphi r
{\si_\infty}$) of $f$ is the usual cokernel (resp. cokernel in $\Maxphi r
{\si_\infty}$) of the inclusion $\ker f \to \frakM$. It follows the
announced property about coimages and then the identification between
image and coimage.
\end{proof}

\begin{lemma}
\label{lem:pmoe}
If $\alpha : \frakM' \to \frakM$ and $\beta : \frakM \to \frakM''$ two
morphisms in $\Maxphi r {\si_\infty}$ such that $\beta \circ \alpha = 0$.
The sequence $0 \to \frakM' \to \frakM \to \frakM'' \to 0$ is exact in
(the abelian category) $\Maxphi r {\si_\infty}$ if and only if the
sequence $0 \to \frakM'[1/u] \to \frakM [1/u] \to \frakM''[1/u] \to 0$
is exact.

Moreover, the functor $\pM_{\O_\E} : \Maxphi r {\si_\infty} \to \Modphi r
{\O_\E}$ is fully faithful.
\end{lemma}

\noindent
{\it Remark:} The reader should be very careful with the following
point. There is two \emph{different} notions of exact sequences in
$\Maxphi r {\si_\infty}$. The first one is given by the structure
of abelian category whereas the second one is just the ``restriction''
of the notion of exact sequence in $\Modphi r {\si_\infty}$. From now on, 
we will only consider the first one. This is for instance the reason why
corollary \ref{cor:maxexact} is not in contradiction with the
counter-example given after proposition \ref{prop:maxleftexact}.

\begin{proof}
By description of kernels and cokernels given in theorem
\ref{theo:abelian}, we have the following: the sequence $0 \to \frakM'
\to \frakM \to \frakM'' \to 0$ is exact in $\Maxphi r {\si_\infty}$ if and
only if $0 \to \frakM' \to \frakM \to \frakM''$ is exact (as a sequence
of $\si$-modules) and $\coker (\frakM \to \frakM'')$ is killed by a
power of $u$. The first part of lemma then follows.

Since for all $\frakM \in \Maxphi r {\si_\infty}$, we have $\frakM \subset
\frakM [1/u]$, the functor $\pM_{\O_\E}$ is clearly faithful. Let
$\frakM$ and $\frakM'$ be two objects of $\Maxphi r {\si_\infty}$ and $f :
\frakM[1/u] \to \frakM'[1/u]$. We have to show that $f$ sends $\frakM$
to $\frakM'$. Using proposition \ref{prop:basickisin} (iii), we have
$f(\frakM) \in \Modphi r {\si_\infty}$ and by the proof of proposition
\ref{prop:supinf}, $f(\frakM) + \frakM'$ (computed in $\frakM'[1/u]$) is 
also an object of $\Modphi r {\si_\infty}$. Hence, by definition of minimal
objects $f(\frakM) + \frakM' \subset \frakM'$, and then $f(\frakM) 
\subset \frakM'$ as required.
\end{proof}

\begin{cor}
\label{cor:fullyfaith}
The functor $T_{\si_\infty}$ defined on $\Maxphi r {\si_\infty}$ is exact
and fully faithful.
\end{cor}

\begin{cor}
\label{cor:maxexact}
The functor $\Max^r : \Modphi r {\si_\infty} \to \Maxphi r {\si_\infty}$ is
exact.
\end{cor}

\begin{theo}
The functor $\Max^r : \Modphi r {\si_\infty} \to \Maxphi r {\si_\infty}$
realizes the localization of $\Modphi r {\si_\infty}$ with respect to
morphisms $f$ such that $T_{\si_\infty}(f)$ is an isomorphism.
\end{theo}

\begin{proof}
Take $\calC$ a category and $F : \Modphi r {\si_\infty} \to \calC$ a
functor that satisfies the following implication: if $T_{\si_\infty}
(f)$ is an isomorphism, then $F(f)$ too. We have to show that there
exists a unique functor $G$ making the following diagram commutative:
$$\xymatrix @C=10pt @R=8pt {
\Modphi r {\si_\infty} \ar[rr]^-{F} \ar[rd]_-{\Max^r} & & \calC \\
& \Maxphi r {\si_\infty} \ar@{.>}[ur]_-{G} }$$
If $\frakM$ is in $\Maxphi r {\si_\infty}$, we must have $G(\frakM) = F
\circ \Max^r(\frakM) = F(\frakM)$. This proves the unicity and gives a
candidate for $G$. Finally, we only have to check that for all $\frakM
\in \Modphi r {\si_\infty}$, there exists a canonical isomorphism between
$F(\frakM)$ and $G(\Max^r(\frakM)) = F(\Max^r(\frakM))$. It is given by
$F(\imax^\frakM)$.
\end{proof}

\subsubsection*{How to recognize maximal objects?}

It seems to be difficult to find a criteria to recognize maximal objects
among objects of $\Modphi r {\si_\infty}$. Nevertheless, we have the
following property of stability.

\begin{prop}
The category $\Maxphi r {\si_\infty}$ is stable under extensions in
$\Modphi r {\si_\infty}$.
\end{prop}

\noindent
{\it Remark.} The proposition means that if $0 \to \frakM' \to \frakM
\to \frakM'' \to 0$ is an exact sequence in $\Modphi r {\si_\infty}$ (and
not in $\Maxphi r {\si_\infty}$ --- that does not make sense) and if
$\frakM'$ and $\frakM''$ are maximal, then $\frakM$ is also. Hence, the
proposition does \emph{not} imply that $\Maxphi r {\si_\infty}$ is the
smallest full subcategory of $\Modphi r {\si_\infty}$ containing simple
objects described in \S \ref{subsec:simple}.

\medskip

\begin{proof}
Assume that $0 \to \frakM' \to \frakM \to \frakM'' \to 0$ is an exact
sequence in $\Modphi r {\si_\infty}$ and $\frakM'$ and $\frakM''$ are
maximal. We have the following diagram:
$$\xymatrix {
0 \ar[r] & \frakM' \ar[r] \ar@{=}[d] & \frakM \ar[r] \ar@{^(->}[d] &
\frakM'' \ar[r] \ar[d]^-{f} & 0 \\
0 \ar[r] & \frakM' \ar[r] & \Max^r(\frakM) \ar[r] & \frakC \ar[r] & 0 }$$
where $\frakC$ is defined as the cokernel of $\frakM' \to
\Max^r(\frakM)$.
A diagram chase shows that $f$ is injective. Moreover by theorem
\ref{theo:abelian}, $\frakC \in \Modphi r {\si_\infty}$ and it is easy to
check that $\frakM''[1/u] = \frakC[1/u]$. Since $\frakM''$ is maximal,
we must have $\frakM'' = \frakC$, \emph{i.e.} $f$ bijective. It follows
that $\frakM = \Max^r(\frakM)$ as required.
\end{proof}

\medskip

Then, we have a sufficient condition to be maximal.

\begin{lemma}
Let $\frakM \in \Modphi r {\si_1}$. If $\coker (\id \otimes \phi)$ is killed by
$u^{p-2}$ then $\frakM$ is maximal.
\end{lemma}

\begin{proof}
It follows from the proof of lemma \ref{lem:boundmax}.
\end{proof}

\begin{cor}
If $er < p-1$, then $\Maxphi r {\si_\infty} = \Modphi r {\si_\infty}$.
\end{cor}

\subsection{Minimal objects}

We develop in this subsection a dual notion of maximal objects (called
\emph{minimal objects}), that satisfies analogous properties. According
to corollary \ref{cor:frsi}, we need to assume $r < \infty$.

\subsubsection*{The functor $\Min^r$}

\begin{deftn}
Let $\frakM \in \Modphi r {\si_\infty}$. The object $\Min^r(\frakM)$ is
defined as the smallest element of $F^r_\si(\frakM[1/u])$. It is endowed
with an homomorphism  $\imin^\frakM : \Min^r(\frakM) \to \frakM$ in the
category $\Modphi r {\si_\infty}$.

An object $\frakM$ of $\Modphi r {\si_\infty}$ is said
\emph{minimal} (in $\Modphi r {\si_\infty}$) if $\imin^\frakM$ is an
isomorphism.
\end{deftn}

\begin{prop}
\label{prop:minfunctor}
The previous definition gives rise to a functor $\Min^r : \Modphi r
{\si_\infty} \to \Modphi r {\si_\infty}$. Moreover, the collection of map
$(\imin^\calM)$ defines a natural transformation between $\Min^r$ and
the identity functor.
\end{prop}

\begin{proof}
Consider $f : \frakM_1 \to \frakM_2$ a map in $\Modphi r {\si_\infty}$. In
order to prove that $\Min^r$ is a functor, we have to show that $f
(\Min^r(\frakM_1)) \subset \Min^r(\frakM_2)$. Since $\Modphi r
{\si_\infty}$ is stable under images (proposition \ref{prop:basickisin}
(iii)), we can assume successively that $f$ is surjective, then
injective.

Assume $f$ surjective. Put $F = f \otimes_{\si} \si[1/u]$ and $\frakM'_1
= F^{-1}(\Min^r(\frakM_2))$. From the surjectivity of $f$ and
$(\Min^r \frakM_2)[1/u] = \frakM_2[1/u]$, we deduce $\frakM'_1[1/u] =
\frakM_1[1/u]$. Moreover, if $\frakK = \ker f$, we have the following
commutative diagram:
$$\xymatrix {
0 \ar[r] & \phi^\star \frakK[1/u] \ar[r] \ar[d]_-{\id \otimes \phi}^-{\sim} &
\phi^\star \frakM'_1 \ar[r] \ar[d]_-{\id \otimes \phi'_1} &
\phi^\star \Min^r(\frakM_2) \ar[r] \ar[d]_-{\id \otimes \phi_2} & 0 \\
0 \ar[r] & \frakK[1/u] \ar[r] & \frakM'_1 \ar[r] & \Min^r(\frakM_2)
\ar[r]  & 0 }$$
Hence $\coker (\id \otimes \phi'_1)$ can be seen as a submodule $\coker
(\id \otimes \phi'_2)$ and so it is killed by $E(u)^r$ (if $r =
\infty$, it must be replaced by a sufficiently large integer).
Therefore, by proposition \ref{prop:basickisin} (ii), $\frakM'_1 \in
F^r_\si(\frakM_1[1/u])$ and $\Min^r(\frakM_1) \subset \frakM'_1$. The
conclusion follows.

Now, assume $f$ injective: we will consider $\frakM_1$ as a subobject of
$\frakM_2$. Put $\frakM'_1 = \frakM_1[1/u] \cap \Min^r(\frakM_2)$. Since 
$(\Min^r \frakM_2)[1/u] = \frakM_2[1/u]$, we have $\frakM'_1[1/u] =
\frakM_1[1/u]$. Now, let
$x \in \frakM'_1$. There exists $y \in \phi^\star \frakM_1[1/u]$ and $z
\in \phi^\star \Min^r(\frakM_2)$ such that $x = \id \otimes \phi(y) =
\id \otimes \phi (z)$. Since $\id \otimes \phi$ is injective on
$\frakM_2[1/u]$, we must have $y = z \in \frakM'_1$. So, by proposition
\ref{prop:basickisin} (ii), $\frakM'_1 \in F^r_\si(\frakM_1[1/u])$.
Hence $\Min^r(\frakM_1) \subset \frakM'_1$, and we are done.

The last statement of the proposition is then obvious.
\end{proof}

\begin{prop}
The functor $\Min^r$ is a projection, that is $\Min^r \circ \Min^r =
\Min^r$.
\end{prop}

\begin{proof}
Just use $\Min^r(\frakM)[1/u] = \frakM[1/u]$.
\end{proof}

\begin{lemma}
\label{lem:imagemin}
Let $f : \frakM \to \frakM'$ a morphism in $\Modphi r {\si_\infty}$. Then
$f(\Min(\frakM)) = \Min(f(\frakM))$.
\end{lemma}

\begin{proof}
First note that $f(\frakM)$ is an object of $\Modphi r {\si_\infty}$
(proposition \ref{prop:basickisin} (iii)) and consequently the formula
$\Min(f(\frakM))$ makes sense.

The inclusion $\subset$ has been proved in proposition
\ref{prop:minfunctor}. Put $\frakM'' = f(\Min(\frakM))$. By proposition
\ref{prop:basickisin} (iii), it is an object of $\Modphi r {\si_\infty}$
such that $\frakM''[1/u] = f(\frakM)[1/u]$. Hence $\Min(f(\frakM))
\subset \frakM''$ as required.
\end{proof}

\begin{cor}
\label{cor:minexact}
Let $f : \frakM \to \frakM'$ a morphism in $\Modphi r {\si_\infty}$. If
$f$ is injective (resp. surjective), then $\Min(f)$ is also.
\end{cor}

\noindent
{\it Remark.} Dualizing the example given after proposition
\ref{prop:maxleftexact}, we see that $\Min$ is not ``middle-exact''.

\begin{prop}
\label{prop:universalmin}
Let $\frakM \in \Modphi r {\si_\infty}$. The couple $(\Min(\frakM),
\imin^\frakM)$ is characterized by the following universal property:
\begin{itemize}
\item the morphism $T_{\si_\infty} (\imin^\frakM)$ is an isomorphism;
\item for each couple $(\frakM', f)$ where $\frakM' \in \Modphi r
{\si_\infty}$ and $f : \frakM' \to \frakM$ becomes an isomorphism under
$T_{\si_\infty}$, there exists a unique map $g : \Min(\frakM) \to
\frakM'$ such that $f \circ g = \imax^\frakM$.
\end{itemize}
\end{prop}

\begin{proof}
The first point is clear. Take $(\frakM', f)$ as in the proposition.
Since $T_{\si_\infty}(f)$ is an isomorphism, $f$ induces an isomorphism
$\frakM'[1/u] \to \frakM[1/u]$ (by full faithfulness of $\pT_{\O_\E}$).
Hence, $f$ is injective, and we can consider $\frakM'$ as a subobject of
$\frakM$. It is then sufficient to prove that $\Min^r(\frakM) \subset
\frakM'$ but this follows from the definition of $\Min^r$.
\end{proof}

\subsubsection*{The category $\Minphi r {\si_\infty}$}

\begin{deftn}
We put $\Minphi r {\si_\infty} = \Min^r(\Modphi r {\si_\infty})$. It is a
full subcategory of $\Modphi r {\si_\infty}$.
\end{deftn}

\begin{prop}
\label{prop:adjonctionmin}
The functor $\Min^r : \Modphi r {\si_\infty} \to \Minphi r {\si_\infty}$ is
a right adjoint of the inclusion functor $\Minphi r {\si_\infty} \to
\Modphi r {\si_\infty}$.
\end{prop}

\begin{proof}
We have to prove that if $f : \frakM \to \frakM'$ is any morphism in
$\Modphi r {\si_\infty}$ with $\frakM$ minimal, then $f$ factors through
$\imin^{\frakM'}$. This is a a direct consequence of proposition
\ref{prop:minfunctor}.
\end{proof}

\begin{theo}
The category $\Minphi r {\si_\infty}$ is abelian. More precisely, if $f :
\frakM \to \frakM'$ is a morphism in $\Minphi r {\si_\infty}$
\begin{itemize}
\item the kernel of $f$ in the usual sense is an object of $\Modphi r
{\si_\infty}$ whose image under $\Min^r$ is a kernel of $f$ in the
abelian category $\Minphi r {\si_\infty}$
\item the cokernel of $f$ in the usual sense, $\coker f$, may have
$u$-torsion; however $\frac{\coker f}{u\text{-torsion}}$ is an object of
$\Minphi r {\si_\infty}$ which is a cokernel of $f$ in the abelian
category $\Minphi r {\si_\infty}$
\item the image (resp. coimage) of $f$ in the usual sense is an object
of $\Minphi r {\si_\infty}$ and is the image (resp. coimage) of $f$ in the
abelian category $\Minphi r {\si_\infty}$.
\end{itemize}
\end{theo}

\begin{proof}
During the proof, we will denote by $\ker f$, $\coker f$, $\im f$ and
$\coim f$ the objects computed in the usual sense.

The assertion about kernels results from propositions
\ref{prop:basickisin} (iii) and \ref{prop:adjonctionmin}. Let's prove
the assertion about cokernels. Denote by $\frakC$ the quotient of
$\coker f$ by its $u$-torsion. Obviously $\frakC$ have no $u$-torsion.
Moreover, it satisfies condition (\ref{eq:condphi}), it is finitely
generated and it is killed by a power of $p$ (since it is a quotient of
$\frakM'$). Hence, by proposition \ref{prop:basickisin} (ii), $\frakC
\in \Modphi r {\si_\infty}$. Lemma \ref{lem:imagemin} applied to the
surjective morphism $\frakM' \to \frakC$ then shows that $\frakC$ is
minimal.

By definition, the image (in $\Minphi r {\si_\infty})$ of $f$, called
$\frakI$, is the kernel (in $\Minphi r {\si_\infty}$) of $\frakM' \to
\frakC$. Hence $\im f \subset \frakI$ and the quotient $\frakI / \im f$
is killed by a power of $u$. It follows that $\Min^r(\im f) =
\Min^r(\frakI)
= \frakI$. But, by lemma \ref{lem:imagemin}, $\im f$ is already minimal.
Thus $\frakI = \im f$ as required. The argument is quite similar for
coimage (remark that since $\coim f$ is isomorphic to $\im f$, it is
also minimal).
\end{proof}

\begin{lemma}
\label{lem:pmoemin}
If $\alpha : \frakM' \to \frakM$ and $\beta : \frakM' \to \frakM''$ two
morphisms in $\Minphi r {\si_\infty}$ such that $\beta \circ \alpha = 0$.
The sequence $0 \to \frakM' \to \frakM \to \frakM'' \to 0$ is exact (in
the abelian category) $\Minphi r {\si_\infty}$ if and only if the sequence
$0 \to \frakM'[1/u] \to \frakM [1/u] \to \frakM''[1/u] \to 0$ is exact.

Moreover, the functor $\pM_{\O_\E} : \Minphi r {\si_\infty} \to \Modphi r
{\si_\infty}$ is fully faithful.
\end{lemma}

\begin{proof}
The first part of lemma follows from the description of kernels and
cokernels given above.

Since for all $\frakM \in \Minphi r {\si_\infty}$, we have $\frakM \subset
\frakM [1/u]$, the functor is clearly faithful. Let $\frakM$ and
$\frakM'$ two objects of $\Minphi r {\si_\infty}$ and $f : \frakM[1/u] \to
\frakM'[1/u]$. We have to show that $f$ sends $\frakM$ to $\frakM'$. The
proof is the same as in proposition \ref{prop:minfunctor}.
\end{proof}

\begin{cor}
\label{cor:fullyfaithmin}
The functor $T_{\si_\infty}$ defined on $\Minphi r {\si_\infty}$ is exact
and fully faithful.
\end{cor}

\begin{cor}
The functor $\Min^r : \Modphi r {\si_\infty} \to \Minphi r {\si_\infty}$ is
exact.
\end{cor}

\subsubsection*{Link with duality}

\begin{prop}
Assume $r$ finite.
For all $\frakM \in \Modphi r {\si_\infty}$, we have natural isomorphisms
$$\Min^r(\frakM^\vee) \simeq \Max^r(\frakM)^\vee
\quad \text{and} \quad
\Max^r(\frakM^\vee) \simeq \Min^r(\frakM)^\vee.$$
In particular, duality permutes subcategories $\Minphi r {\si_\infty}$
and $\Maxphi r {\si_\infty}$.
\end{prop}

\begin{proof}
Formula (\ref{eq:tsidual}) implies that, given a morphism
$f$ in the category $\Modphi r {\si_\infty}$, $T_{\si_\infty} (f)$ is an
isomorphism if and only if $T_{\si_\infty} (f^\vee)$ is.
Then, the proposition is a formal (and easy) consequence of the
universal properties defining $\Max^r$ (proposition
\ref{prop:universal}) and $\Min^r$ (proposition \ref{prop:universalmin})
on the one hand, and the full faithfulness of $T_{\si_\infty}$ on $\Maxphi 
r {\si_\infty}$ (corollary \ref{cor:fullyfaith}) and $\Minphi r
{\si_\infty}$ (corollary \ref{cor:fullyfaithmin}) on the other hand.
\end{proof}

\subsection{A reciprocity formula}

In this subsection, we will use the functor $j_\star$ of Fontaine
defined in \S B.1.4 of \cite{fontaine-fest}. For $M \in \pModphi {}
{\O_\E}$, define the ordered set $G_\si (M)$ as the set of
$\si$-submodules $\frakM \subset M$ such that $\frakM$ is of finite
type over $\si$, stable under $\phi$ and $\id \otimes \phi :
\phi^\star \frakM[1/u] \to \frakM[1/u]$ is bijective. Recall that,
by definition:
$$j_\star M = \bigcup_{\frakM \in G_\si(M)} \frakM.$$
In the same way, we put for any $r \in \{0, 1, \ldots, \infty\}$:
$$j^r_\star M = \bigcup_{\frakM \in G^r_\si(M)} \frakM$$
where $G^r_\si(M)$ is the ordered set of all $\frakM \in \Modphi r
{\si_\infty}$ with $\frakM \subset M$ (we do not ask $\frakM[1/u]$ to be
equal to $M$). By \S B.1.5.3 of \cite{fontaine-fest}, the equality
$G_\si(M) = G^\infty_\si(M)$ holds.
Moreover, if $\frakM$ is an object of $\Modphi r {\si_\infty}$, (the
proof of) proposition \ref{prop:supinf} shows that greatest elements of
$F^r_\si(M)$ and $G^r_\si(M)$ coincide. Hence $\Max^r(\frakM) =
j_\star^r(\frakM[1/u])$.

Following \cite{liu2}, we define for $r \in \{0, 1, \ldots, \infty\}$:
$$\si_n^{f,r} = j^r_\star (\O_{\E^\ur}/p^n\O_{\E^\ur}) \subset
\O_{\E^\ur}/p^n\O_{\E^\ur} \quad \text{and} \quad
\si^{f,r} = \varprojlim_n \si_n^{f,r} \subset \O_{\E^\ur}.$$
For all integer $n$, $\si_n^{f,r}$ is an object of $\pModphi r
{\si_\infty}$, and obviously $\si_n^{f,\infty} = \bigcup_{r \in \N}
\si_n^{f,r}$. By proposition 2.5.1 of \emph{loc. cit.}, they are stable
under $\phi$ and the action of $G_\infty$. Furthermore, this proposition
implies that $\si^{f,\infty}$ is the period ring $\si^\ur$
traditionally used in this context (for instance in \cite{kisin05},
\cite{liu}, \cite{liu2}). Finally, if $\frakM \in \Modphi r
{\si_\infty}$ is cancelled by $p^n$, the formula for $T_{\si_\infty}
(\frakM)$ can be ``simplified'' as follows:
$$T_{\si_\infty}(\frakM) = \hom_{\pModphi r \si} (\frakM,
\si_n^{f,r}).$$
(To prove this, it is enough to remark that the image of any $f \in
T_{\si_\infty} (\frakM)$ is an object of $\Modphi r {\si_\infty}$, which
follows more or less from proposition \ref{prop:basickisin} (iii).)

\bigskip

Here is the main theorem of this subsection:

\begin{theo}
Let $\frakM \in \Modphi r {\si_\infty}$ killed by $p^n$. Then
$\Max^r(\frakM) = \hom_{\Z_p[G_\infty]} (T_{\si_\infty}(\frakM),
\si_n^{f,r})$.
\end{theo}

\noindent
{\it Remark.} It seems that such a formula does not exist with $\Min^r$
(instead of $\Max^r$). Indeed, it would probably imply the
left-exactness of $\Min^r$, which is known to be false (see remark after
corollary \ref{cor:minexact}).

\begin{proof}
Put $\tilde \frakM = \hom_{\Z_p[G_\infty]} (T_{\si_\infty}(\frakM),
\si_n^{f,r})$. It is endowed with a Frobenius $\phi$ (given by the
Frobenius on $\si_n^{f,r}$). Moreover, biduality gives a natural map
compatible with Frobenius:
$$\iota : \Max^r(\frakM) \to \hom_{\Z_p[G_\infty]}
(T_{\si_\infty}(\Max^r(\frakM)), \si_n^{f,r}) \simeq \tilde \frakM.$$
By remark A.1.2.7.(a) of \cite{fontaine-fest}, the composite
$$\xymatrix @C=10pt @R=0pt {
\frakM[1/u] \ar[rr]^-{\iota \otimes_\si \si[1/u]} & \hspace{1cm}
& \tilde \frakM[1/u] \ar@{^(->}[r] & \hom_{\Z_p[G_\infty]}
(T_{\si_\infty}(\Max^r(\frakM)), \O_{\E^\ur} / p^n \O_{\E^\ur}) \\
& & & \hspace{1cm} = \hom_{\Z_p[G_\infty]}
(\pT_{\O_\E}(\frakM[1/u]), \O_{\E^\ur} / p^n \O_{\E^\ur}) }$$
is bijective. Hence, $\iota \otimes_\si \si[1/u]$ is also a bijection.
We want to prove that $\iota$ itself is an isomorphism. Injectivity is
clear since $\Max^r(\frakM)$ have no $u$-torsion. Since $\Max^r(\frakM)
= j^r_\star(\frakM[1/u])$, surjectivity will follow from the statement
``every $f \in \tilde \frakM$ is contained in an
object $\frakN \in G^r_\si(\frakM[1/u])$''. Let us prove the claim.
Consider $e_1, \ldots, e_d$ a generating family of $\frakM$ and put $x_i
= f(e_i)$. By definition of $\si_n^{f,r}$, there exists $\frakN_i
\subset \O_{\E^\ur}/p^n \O_{\E^\ur}$ with $\frakN_i \in \Modphi r
{\si_\infty}$ and $x_i \in \frakN_i$. Then, as usual using proposition
\ref{prop:basickisin}, we can check that $\frakN =
\hom_{\Z_p[G_\infty]} (T_{\si_\infty}(\Max^r(\frakM)), \sum_{i=1}^d
\frakN_i)$ answers the question.
\end{proof}

\begin{cor}
If $\frakM$ a simple object of the abelian category $\Maxphi r
{\si_\infty}$, then $T_{\si_\infty} (\frakM)$ is an irreducible
representation.
\end{cor}

\begin{cor}
For any $r$, the (essential closure of the) category
$T_{\si_\infty}(\Modphi r {\si_\infty})$ is stable under quotients and
subobjects.
\end{cor}

\begin{proof}
Noting that $T_{\si_\infty}(\Modphi r {\si_\infty}) =
T_{\si_\infty}(\Maxphi r {\si_\infty})$, the corollary is a direct
consequence of property 6.4.2 of \cite{caruso-crelle}.
\end{proof}

\subsection{Simple objects}
\label{subsec:simple}

For simplicity, we assume in this subsection $c_0 = 1$ (recall that $c_0
= \frac{E(0)} p)$. Of course, it is not crucial but assuming this will
allow us to simplify several formulas and several definitions of
objects.

We fix an element $r \in \{0, 1, 2, \ldots, \infty\}$.

\paragraph{Definitions and basic properties}

\begin{deftn}
\label{def:simple}
Let $\calS'$ be the set of sequences of integers between $0$ and $er$
that are periodic (from the start). To a sequence $(n_i) \in \calS$, we 
associate several numeric invariants:
\begin{itemize}
\item its dimension $d$: it is the smallest period of $(n_i)$;
\item for $i \in \Z/d\Z$, the integer $s_i = n_i p^{d-1} + n_{i+1}
p^{d-2} + \cdots + n_{i+d-1}$;
\item for $i \in \Z/d\Z$, $t_i = \frac{s_i}{p^d -1} \in \Q/\Z$ and $t =
t_0$.
\end{itemize}
We also associate an object $\frakM(n_i) \in \Modphi r {\si_1}$ defined
as follows:
\begin{itemize}
\item as a $\si_1$-module, $\frakM(n_i) = \displaystyle \bigoplus_{i \in
\Z/d\Z} e_i \: \si_1$;
\item for all $i \in \Z/d\Z$, $\phi(e_i) = u^{n_i} e_{i+1}$.
\end{itemize}
Let $\calS$ be the subset of $\calS'$ consisting of all sequences $(n_i)$
for which the elements $t_0, \ldots, t_{d-1}$ are pairwise distinct (in
$\Q/\Z$).
\end{deftn}

\begin{prop}
Assume $r < \infty$.
Let $(n_i)$ and $(m_i)$ be two sequences in $\calS'$. If $n_i + m_i =
er$ for all $i$, then duality permutes objects $\frakM(n_i)$ and
$\frakM(m_i)$.
\end{prop}

\begin{proof}
Easy computation.
\end{proof}

\begin{lemma}
\label{lem:eqphi}
Let $(n_i) \in \calS$ and $s$ be a non negative integer. Let $(E)$ be
the equation $\phi^d(x) = u^s x$ in variable $x \in \frakM(n_i)$ (resp.
$x \in \frakM(n_i)[1/u]$). Then $(E)$ has a non zero solution if and
only if there exists $i \in \Z/d\Z$ (necessary unique) and $v$ a non
negative integer (resp. an integer) such that $s-s_i = v(p^d-1)$. In
this case, the set of solutions is $\{\alpha u^v e_i, \: \alpha \in k
\cap \F_{p^d}\}$.
\end{lemma}

\begin{proof}
First, remark that if $p^d-1$ divides $s - s_i$ and $s - s_j$, we get
$s_i \equiv s_j \pmod {p^d-1}$ and then $t_i \equiv t_j \pmod \Z$.
Hence, by assumption, $i = j$ (in $\Z/d\Z$). This justifies the unicity
of $i$.

An easy computation gives $\phi^d (e_i) = u^{s_i} e_i$ for all $i$.
Write $x = x_0 e_0 + \cdots + x_{d-1} e_{d-1}$ with $x_i \in \si_1 =
k[[u]]$ (resp $x_i \in \si[1/u]$). Then, the equation $(E)$ becomes the
system $u^{s_i} x_i^{p^d} = u^s x_i$, and the lemma follows.
\end{proof}

\begin{prop}
\label{prop:isomsimple}
Let $(n_i)$ and $(n'_i)$ be in $\calS$. The objects $\frakM(n_i)$ and
$\frakM(n'_i)$ are isomorphic if and only if there exists an integer $b$
such that $n'_{i+b} = n_i$ for all $i$.
\end{prop}

\begin{proof}
The condition is obviously sufficient. Now, take
$(n_i)$, $d$ and $s_i$, \emph{etc}. as in the definition \ref{def:simple}. We
have to show that knowing $\frakM = \frakM(n_i)$, we can recover the
sequence $(n_i)$ up to a shift. Since $d$ is the dimension of $\frakM$,
it is clearly determined. Remark that by lemma \ref{lem:eqphi}, 
integers $s_i$ are exactly integers $s$ for which there exists $x
\in \frakM$, $x \not\in u \frakM$ such that $\phi^d(x) = u^s x$. So,
their set is also determined. Moreover if $x_i$ is a non zero solution
of $\phi^d(x_i) = u^{s_i} x_i$, we can write $x_i = \alpha_i e_i$ with
$\alpha_i \in k$. It follows that $\phi$ maps $x_i \si_1$ to
$x_{i+1} \si_1$ and then that the sequence $(s_i)$ is determined up to
circular permutation. It remains to prove that the knowledge of $(s_i)$
determines the sequence $(n_i)$. But we have an equality
$$\pa{ \begin{array}{c} s_0 \\ s_1 \\ \vdots \\ s_{d-1} \end{array} } =
M \pa{ \begin{array}{c} n_{d-1} \\ n_0 \\ \vdots \\ n_{d-2} \end{array}
}$$
where $M$ is a matrix with integer coefficients whose reduction modulo
$p$ is identity. The proposition follows.
\end{proof}

\paragraph{Maximum and minimum objects}

Here, we compute functors $\Min^r$ and $\Max^r$ on objects
$\frakM(n_i)$. We first define several subsets of $\calS'$.

\begin{deftn}
\label{def:minmaxS}
Put $m = \min\{er, p-1\}$.

Let $\calS_\max \subset \calS'$ be the set of sequences of integers
between $0$ and $m$ that are periodic except that the constant sequence
with value $p-1$ is removed from $\calS_\max$ (if necessary).

If $r < \infty$, define $\calS_\min \subset \calS'$ as the set of
sequences of integers between $er-m$ and $er$ that are periodic except
that the constant sequence with value $er-(p-1)$ is removed from
$\calS_\min$ (if necessary).
\end{deftn}

\begin{lemma}
We have $\calS_\max \subset \calS$ and $\calS_\min \subset \calS$ (if
$r$ is finite).
\end{lemma}

\begin{proof}
Exercise. (For $\Max$, one may consider expansion of $t_i$'s in 
$p$-basis.)
\end{proof}

Until the end of this subsection, the assumption $r < \infty$ will
always be implicit when dealing with minimal objects.

\begin{prop}
\label{prop:objmax}
Let $(n_i) \in \calS_\max$ (resp $(n_i) \in \calS_\min)$. Then,
$\frakM(n_i)$ is maximal (resp. minimal).
\end{prop}

\begin{proof}
By duality, we only have to prove the statement with $\Max$.
By examining the proof of lemma \ref{lem:boundmax}, we see that
$\Max(\frakM(n_i)) \subset
\frac 1 u \frakM(n_i)$. Assume by contradiction, that there exists an
element $x \in \Max(\frakM(n_i))$, $x \not\in \frakM(n_i)$ and write $ux
= x_0 e_0 + \cdots + x_{d-1} e_{d-1}$ with $x_i \in \si_1$ and $x_j
\not\in u \si_1$ for one index $j$. A computation gives:
$$\phi(x) = \frac{\phi(x_0)}{u^{p-n_0}} e_1 + \cdots +
\frac{\phi(x_{d-2})}{u^{p-n_{d-2}}} e_{d-1} +
\frac{\phi(x_{d-1})}{u^{p-n_{d-1}}} e_0.$$
This element have to lie in $\Max(\frakM(n_i))$, which implies $p-n_j
\leq 1$, \emph{i.e.} $n_j \geq p-1$. So $n_j = p-1$. Repeating the
argument with $\phi(x)$ instead of $x$, we obtain $n_{j+1} = p-1$, and
so on. Finally, $n_i = p-1$ for all $i$ and $(n_i) \not\in \calS_\max$.
\end{proof}

\begin{prop}
\label{prop:compmax}
For any $(n_i) \in \calS$, there exists a sequence $(m_i) \in
\calS_\max$ (resp. $(m_i) \in \calS_\min$) such that $\Max(\frakM(n_i))
= \frakM(m_i)$ (resp. $\Min(\frakM(n_i)) = \frakM(m_i)$).
\end{prop}

\begin{proof}
By duality, we only have to prove the statement with $\Max$. Denote by
$s'_i$ the unique integer in $[0, p^d-1[$ congruent to $s_i$ modulo
$p^d-1$, and define $m_i$ to be the quotient in the Euclidean division of
$s'_i$ by $p$. It is easy to see that the $m_i$'s ($0 \leq i \leq
d-1$) are digits in
$p$-basis of $s'_0$, and that this property implies $(m_i) \in
\calS_\max$. Now, put $q_i = \frac{s_i - s'_i}{p^d-1}$: it is the
quotient in the Euclidean division of $s_i$ by $p$. These numbers
are non negative integers and they satisfy the relation $p q_i + m_i =
q_{i+1} + n_i$ for all $i \in \Z/d\Z$.

Denote by $\frakM'$ the submodule of $\frakM[1/u]$ generated by $e'_i =
\frac 1{u^{q_i}} e_i$. A direct computation gives $\phi(e'_i) = u^{m_i}
e'_{i+1}$, and then $\frakM' \simeq \frakM(m_i)$. Moreover proposition
\ref{prop:objmax} shows that $\frakM'$ is maximal. The conclusion
follows.
\end{proof}

\noindent
{\it Remark.}
If $(n_i)$ is in $\calS'$ but not in $\calS$, almost all arguments of
the proof are still correct. The only problem is that the sequence
$(m_i)$ obtained is periodic with period less than $d$.

\begin{cor}
Let $(n_i) \in \calS$. If $\frakM(n_i)$ is maximal (resp. minimal) then,
$(n_i)$ is in $\calS_\max$ (resp. $\calS_\min$).
\end{cor}

\begin{proof}
By proposition \ref{prop:compmax}, we can find a sequence $(m_i) \in
\calS_\max$ such that $\frakM(n_i) = \Max(\frakM(n_i)) \simeq
\frakM(m_i)$. By proposition \ref{prop:isomsimple}, there exists an
integer $b$ such that $n_i = m_{i+b}$ for all $i$, and then $(n_i) \in
\calS_\max$.
\end{proof}

\begin{cor}
\label{cor:maxisom}
Let $(n_i)$ and $(n'_i)$ be in $\calS$. Objects $\Max(\frakM(n_i))$
(resp. $\Min(\frakM(n_i))$) and $\Max(\frakM(n'_i))$ (resp.
$\Min(\frakM(n'_i))$) are isomorphic if and only if there exists an
integer $b$ such that $t \equiv p^b t' \pmod \Z$ (with obvious
notations).
\end{cor}

\begin{proof}
Easy after proposition \ref{prop:isomsimple} and proof of proposition
\ref{prop:compmax}.
\end{proof}

\paragraph{Classification}

With notations of \S 1 of \cite{serre}, an easy computation gives the
following theorem.

\begin{theo}
We assume $k$ to be algebraically closed. Let $(n_i) \in \calS_\max$.
Then $T_{\si_\infty}(\frakM(n_i))$ is an irreducible representation of
$G_\infty$ whose tame inertia weights are exactly the $n_i$'s.
\end{theo}

\noindent
{\it Remark.} For $(n_i) \in \calS_\min$, tame inertia weights of
$T_{\si_\infty}(\frakM(n_i))$ are not simply linked with the $n_i$'s.
Precisely, to make the computation, the method is to write the rational
number $t_i$ in $p$-basis and then to read its digits.

\begin{prop}
\label{prop:simple}
We assume $k$ to be algebraically closed.
Let $(n_i) \in \calS$. The object $\Max(\frakM(n_i))$ (resp.
$\Min(\frakM(n_i))$) is simple in $\Maxphi r {\si_\infty}$ (resp.
$\Minphi r {\si_\infty}$). All simple objects can be written in this
form.
\end{prop}

\begin{proof}
If $er < p-1$, the proposition was already proved in \S 4 of
\cite{caruso-crelle}. From now on, we assume $er \geq p-1$. Moreover, it
suffices, using duality, to show the proposition with $\Max$.

By the exactness and the full faithfulness of $T_{\si_\infty}$ on
$\Maxphi r {\si_\infty}$ (corollary \ref{cor:fullyfaith}), in order to
show that $\Max(\frakM(n_i))$ is simple, it is enough to justify that
$T_{\si_\infty}(\Max(\frakM(n_i)))$ is an irreducible representation,
which is a direct consequence of the previous theorem. Now, consider
$\frakM \in \Maxphi r {\si_\infty}$ a simple object. By the previous
theorem and the classification of irreducible representations given in
\S 1.5 and \S 1.6 of \cite{serre}\footnote{In this reference, the
classification is made for $G_K$-representations, but it is easily seen
that the same arguments works with $G_\infty$-representations.}, there
exists a quotient of $T_{\si_\infty} (\frakM)$ isomorphic to
$T_{\si_\infty}(\frakM(n_i))$ for some sequence $(n_i) \in
\calS_\max$. Since $er \geq p-1$, we have $\frakM(n_i) \in \Modphi r
{\si_\infty}$ and $\frakM(n_i) = \Max^r(\frakM(n_i))$ (since $(n_i)$ is in
$\calS_\max$). Finally, full faithfulness of $T_{\si_\infty}$ on
$\Maxphi r {\si_\infty}$ gives a non-vanishing morphism $\frakM(n_i) 
\to \frakM$, and the proposition follows.

\medskip

Instead of using properties of $T_{\si_\infty}$, we can translate
Serre's proof to obtain a classification of simple objects of $\pModphi
{} {\O_\E}$ (which then implies easily the proposition). Since it seems
difficult to find a reference for this classification, we give it here.

Let $M$ be a simple object in $\pModphi {} {\O_\E}$. We will prove that $M$
is isomorphic to $\frakM(n_i)[1/u]$ for a sequence $(n_i) \in
\calS_\max$. First remark that simplicity shows directly that $M$ is killed
by $p$, and hence is $k((u))$-vector space. Let's call $\calL(M)$ the
$k((u))$-vector space of all $k((u))$-linear endomorphisms of $M$ and
denote by $E$ the subset of $\calL(M)$ consisting of those that commute
with Frobenius. Since $M$ is simple, Schur lemma implies that $E$ is a
field. Moreover, it is an $\F_p$-vector space and we have a canonical
$k((u))$-linear map $\alpha : k((u)) \otimes_{\F_p} E \to \calL(M)$. We
claim that $\alpha$ is injective. Indeed, consider $(f_i)_{i \in
I}$ a basis (not necessarly finite) of $E$ over $\F_p$ and assume by
contradiction that $\ker \alpha \neq 0$. Consider an element $f \in \ker
\alpha$ written $f = \sum_{j \in J} a_j \otimes f_j$ where $J \subset I$
is finite and not empty, and where $a_j \neq 0$ for all $j \in J$. Assume 
moreover that $\card J$ is minimal.
Applying Frobenius to $f$, we find $f^\phi = \sum_{j \in J} a_j^p \otimes 
f_j \in \ker \alpha$. Since $\alpha_{|E}$ is obviously injective, it is
impossible that all the $a_j$'s are congruent modulo $\F_p^\star$. Hence, 
a suitable linear
combination of $f$ and $f^\phi$ gives a non-trivial element in $\ker
\alpha$ that can be written $\sum_{j \in J'} b_j \otimes f_j$ with
$J' \subsetneq J$, $J \neq \emptyset$, contradicting the minimality of 
$\card J$ and proving the claim.

It follows that $E$ is finite dimensional over $\F_p$ and then
himself finite. Thus, $E$ is a finite field. In particular, by
Wedderburn's theorem, it is commutative. Moreover, by definition, it
acts on $M$, making $M$ a module over $E \otimes_{F_p} k((u))$. Since
$k$ is algebraically closed, this tensor product splits completely.
Precisely, if $d$ is the degree of $E$ over $\F_p$, we have an 
isomorphism $E \otimes_{\F_p} k((u)) \simeq k((u))^d$, $x \otimes y
\mapsto (x^{p^{-i}} y)_{i \in \Z/d\Z}$. Considering idempotents of this
decomposition, we have a canonical splitting $M = M_1 \oplus M_2
\oplus \cdots \oplus M_d$ where $M_i$ is a vector space over $k((u))$.
Examining the semi-linearity of $\phi$, it is easily seen that $\phi$ 
maps $M_i$ to $M_{i+1}$. Consequently $\phi^d$ maps $M_1$ to himself,
and since $k$ is algebraically close, it must exist an eigenvector $E_1$
of $\phi^d : M_1 \to M_1$, say $\phi^d (E_1) = \lambda E_1$ with 
$\lambda \neq 0$ by étaleness of $M$. Replacing
$E_1$ by $\mu E_1$ changes $\lambda$ into $\mu^{p^d-1} \lambda$. This
allows us to assume that $\lambda = u^s$ for an integer $s \in \{0, 1,
\ldots, p^d-2\}$. Writing $s$ in $p$-basis, we have $s = n_1 p^{d-1}
+ n_2 p^{d-2} + \cdots + n_d$ for some sequence $(n_i) \in \calS_\max$.
Now, we define further $E_i$'s by the inductive
formula $E_{i+1} = u^{-n_i} \phi(E_i)$. A simple computation gives
$E_{d+1} = E_1$. Finally, if $d'$ is the smallest period of $(n_i)$ (which is
a divisor of $d$), it remains easy to check that the map $\frakM(n_i)[1/u] 
\to M$, $e_i \mapsto E_i + E_{i+d'} + E_{i + 2d'} \cdots + E_{i+d-d'}$
is an injective morphism in $\pModphi {} {\O_\E}$. Since $M$ is simple, it 
is an isomorphism and we are done.
\end{proof}

\subsection[Reformulation with $\text{\rm
Mod}^{r,\phi,N}_{S_\infty}$]{Reformulation with $\Modphi r
{S_\infty}$}

Under the equivalence of the theorem \ref{theo:equivtors}, previous
results imply theorem \ref{theo:overS} of the introduction. Moreover,
with notations of theorem \ref{theo:overS}, duality on $\Modphi r
{S_\infty}$ discussed in \S \ref{subsec:duality} permutes functors
$\Max^r$ and $\Min^r$ and categories $\Maxphi r {S_\infty}$ and $\Minphi
r {S_\infty}$ (here $r < p-1$).

\medskip

Furthermore, if $k$ is algebraically close, we have a classification of
simple objects of $\Maxphi r {S_\infty}$ and $\Minphi r {S_\infty}$. For
any sequence $(n_i) \in \calS$ (see definitions \ref{def:simple}) put
$\calM(n_i) = M_{\si_\infty} (\frakM(n_i))$. It is described as follows:
\begin{itemize}
\item $\calM(n_i) = \displaystyle \bigoplus_{i \in \Z/d\Z} f_i \: S_1$;
\item $\Fil^r \calM(n_i) = \displaystyle \sum_{i \in \Z/d\Z} u^{er-n_i}
f_i \: S_1$;
\item for all $i \in \Z/d\Z$, $\phi_r(u^{er-n_i} f_i) = (-1)^r
f_{i+1}$.
\end{itemize}

\begin{theo}
\label{theo:simple}
Assume the residue field $k$ algebraically closed, and $r < p-1$.

For all sequence $(n_i) \in \calS_\max$ (resp. $(n_i) \in \calS_\min$),
the object $\calM(n_i)$ is simple in $\Maxphi r {S_\infty}$ (resp. in
$\Minphi r {S_\infty}$). Every simple object of $\Maxphi r {S_\infty}$
(resp. of $\Minphi r {S_\infty}$) is isomorphic to $\calM(n_i)$ for a
certain sequence $(n_i) \in \calS_\max$ (resp. $(n_i) \in \calS_\min$).
Moreover, two objects $\calM(n_i)$ and $\calM(m_i)$ are isomorphic if
and only if there exists an integer $b$ such that $n_i = m_{i+b}$ for
all $i$.

The $G_\infty$-representation $T_\qcris(\calM(n_i))$ is irreducible and
its tame inertia weights are exactly the $n_i$'s.
\end{theo}

\section{The case $r=1$}

We assume $r=1$.
The forgetting functor $\Modphin 1 {S_\infty} \to \Modphi 1 {S_\infty}$
is an equivalence of categories (see lemma 5.1.2 of \cite{bcdt}), and
therefore, quasi-semi-stable representations are exactly restrictions
to $G_\infty$ of quotients of two lattices in a crystalline
representation with Hodge-Tate weights in $\{0,1\}$. Moreover, they are
also (restrictions to $G_\infty$ of) representations of the form
$\calG(\bar K)$ where $\calG$ is a finite flat group scheme over $\O_K$
killed by a power of $p$. Let denote by $\Repg (G_K)$ (resp. $\Repg
(G_\infty)$) their category. We have the following commutative diagram
$$\xymatrix{
\Modphin 1 {S_\infty} \ar[d]_-{\sim} \ar[rr]^-{T_\st} & & \Repg(G_K)
\ar[d] \\
\Modphi 1 {S_\infty} \ar[r]^-{\Max^1} & \Maxphi 1 {S_\infty}
\ar[r]^-{T_\qcris} & \Repg(G_\infty) }$$
where vertical arrows represent forgetting functors.

\begin{prop}
The functor $T_\st$ factors through $\Maxphi 1 {S_\infty}$.
\end{prop}

\begin{proof}
By the last statement of theorem \ref{theo:overS}, it is sufficient to
prove that if $T_\qcris(f)$ is an isomorphism, then $T_\st (f)$ is also
(where $f$ in any map in $\Maxphi 1 {S_\infty}$). But it is obvious
since $T_\qcris(f) = T_\st(f)$.
\end{proof}

\begin{cor}
The functor $\Repg (G_K) \to \Repg (G_\infty)$ is fully faithful. In
other words, if $F : T \to T'$ is a $G_\infty$-equivariant map between
two objects of $\Repg(G_K)$, then it is $G_K$-equivariant.

Moreover, $T_\st : \Maxphi 1 {S_\infty} \to \Repg(G_K)$ is fully
faithful.
\end{cor}

\begin{proof}
If $\calM$ and $\calM'$ are objects of $\Maxphi 1 {S_\infty}$, the
composite $$\hom_{\Maxphi 1 {S_\infty}} (\calM, \calM') \to
\hom_{\Repg(G_K)} (T_\st(\calM'), T_\st(\calM)) \to
\hom_{\Repg(G_\infty)} (T_\qcris(\calM'), T_\qcris(\calM))$$ is
bijective (by full faithfulness of $T_\qcris$) whereas the second map is
obviously injective. This implies that both maps are bijective. Since
$T_\st : \Maxphi 1 {S_\infty} \to \Repg (G_K)$ is essentially surjective
(by definition of $\Repg (G_K)$), the corollary follows.
\end{proof}

\noindent
{\it Remark.} The first part of corollary was already known (theorem
3.4.3 of \cite{breuil-azumino}). However, the proof given here is
slightly different.

\section{Perspectives and questions}

\subsubsection*{The semi-stable and crystalline case}

Of course, one may ask if the previous theory can be extended to the
semi-stable case. Precisely:

\begin{question}
\label{quest:quot}
Can we find a simple criteria to recognize an object of $\ModphiN r
{S_\infty}$ that can be written as a quotient of two strongly divisible
modules?
\end{question}

\begin{question}
\label{quest:withN}
Are theorems \ref{theo:overS} and \ref{theo:simple} (with $N(f_i) = 0$)
still true if we replace $\Modphi r {S_\infty}$ by $\ModphiN r {S_\infty}$
($\Maxphi r {S_\infty}$ by $\MaxphiN r {S_\infty}$, and $T_\qcris$ by
$T_\st$)?
\end{question}

It seems quite difficult to find a satisfying answer to question
\ref{quest:quot}. For the moment, the authors do not know if any object
can be written such as a quotient, although they conjecture it is false.
On the other hand, question \ref{quest:withN} seems more accessible and
will be partially answered in a forthcoming paper.

\medskip

Finally note that links between crystalline and semi-stable torsion
theory seem to be more complicated than it looks. Denote by $\Modphin r
{S_\infty}$ the full subcategory of $\ModphiN r {S_\infty}$ gathering
objects $\calM$ satisfying $N(\calM) \subset (uS + \Fil^1 S) \calM$. If
$r=1$, we saw that the forgetting functor $\Modphin r {S_\infty} \to
\Modphi r {S_\infty}$ is an equivalence and then allows us to identify
$\Modphin r {S_\infty}$ and $\Modphi r {S_\infty}$. However, if $r > 1$,
this functor is not anymore fully faithful and consequently one can
\emph{not} identify $\Modphin r {S_\infty}$ as a subcategory of $\Modphi r
{S_\infty}$.

Here is a counter-example. Assume $e \geq \frac{p-1}{r-1}$. Assume also
that there exists $\lambda \in S_1$ such that $\lambda^{p-1} \equiv c
\pmod p$. Put $\calM = e_1 S_1 \oplus e_2 S_1$, and let $\Fil^r \calM$
be the
submodule of $\calM$ generated by $e_1$, $u^{e+p-1} e_2$ and $\Fil^p S_1
\calM$. Equip $\calM$ with a Frobenius by putting $\phi_r(e_1) = e_1$
and $\phi_r (u^{e+p-1} e_2) = e_2$. Then, it is possible to define on
$\calM$ two monodromy operators $N_1$ and $N_2$ by the formulas
$N_1(e_1) = N_2(e_1) = 0$, $N_2(e_1) = \lambda u^p e_2$, $N_2(e_2) = 0$.
These operators give rise to two objects $\calM_1$ and $\calM_2$ of
$\Modphin r {S_\infty}$. They are not isomorphic since $N \circ \phi_r$
vanishes on $\Fil^r \calM_1$ but not on $\Fil^r \calM_2$. Moreover, one
can prove that associated Galois representations (\emph{via} the functor
$T_\st$) are not isomorphic.

Going further, we can evaluate what should be $\Min(\calM_1)$ and
$\Min(\calM_2)$. For simplicity, assume $e < p-1$. Define $\calM' =
e'_1 S_1 \oplus e'_2 S_1$ endowed with $\Fil^r \calM'$ generated by
$e'_1$, $u^e e'_2$ and $\Fil^p S \calM'$. Put $\phi_r(e'_1) = e'_1$ and
$\phi_r(u^e e'_2) = e'_2$. Again, we can equip $\calM'$ with two
monodromy operators $N_1$ and $N_2$ defined by $N_1(e'_1) =
N_1(e'_2) = 0$, $N_2(e'_1) = \lambda e'_2$ and $N_2(e'_2) = 0$. Call
$\calM'_1$ and $\calM'_2$ the corresponding objects of $\ModphiN r
{S_\infty}$. For $i \in \{1,2\}$, we have a morphism $\calM'_i \to
\calM_i$ (in $\ModphiN r {S_\infty}$) and we can check that it induces an
isomorphism \emph{via} $T_\st$. Moreover, since $e \leq p-2$, $\calM'_1$
and $\calM'_2$ should be minimal. Therefore $\Min^r(\calM_i)$ should be
equal to $\calM'_i$ and the implication $(\calM \in \Modphin r
{S_\infty}) \Rightarrow (\Min(\calM) \in \Modphin r {S_\infty})$ should
(surprisingly) be false.

\subsubsection*{A point of view with sheaves}

Proposition \ref{prop:adjonction} and theorem \ref{theo:abelian} show
that the situation is quite similar to what happens with presheaves and
sheaves. More concretely we may ask the following question:

\begin{question}
\label{quest:sheaves}
Is it possible to see objects of $\Modphi r {\si_\infty}$ (resp. $\Maxphi r
{\si_\infty}$) as global sections of some presheaves (resp. sheaves) on
a certain site, in such a way that the functor $\Max$ corresponds to the
functor ``associated sheaf''?

Is it possible to find such presheaves and sheaves in certain cohomology
groups of certain varieties?
\end{question}

In order to precise the latest question, assume $r = 1$. Consider
$\calG$ a finite flat group scheme killed by a power of $p$ over $\O_K$.
In \cite{breuil-group}, Breuil manages to associate to $\calG$ an object
$\calM \in \Modphi r {S_\infty}$ using geometric construction. We can ask
the following:

\begin{question}
Is it possible to find an only geometric recipe that associates to
$\calG$ the object $\Max(\calM)$? For instance, can we obtain this
recipe by sheafifying (in a certain way) the construction of Breuil?
\end{question}

\end{document}